\documentclass[a4paper, 11pt]{article}
\usepackage{bbm}
\usepackage{amsfonts}
\usepackage{amsmath}
\usepackage{amsthm,mathrsfs}
\usepackage{amsfonts,amssymb}
\usepackage{enumerate}
\usepackage{xcolor}
\usepackage{setspace}
\usepackage{authblk,titletoc}
\usepackage{xeCJK}
\allowdisplaybreaks[4]
\newtheorem{theorem}{Theorem}[section]
\newtheorem{proposition}[theorem]{Proposition}%[section]

\newtheorem{lemma}[theorem]{Lemma}%[section]
\newtheorem{remark}[theorem]{Remark}%[section]

\numberwithin{equation}{section}

\def\pat{\partial}

\def\az{\alpha}
\def\be{\beta}

\def\A{\mathcal{A}}

\def\O{\mathcal{O}}

\def\ep{\epsilon}

\def\pat{\partial}
\def\De{\Delta}

\def\rr{{\mathbb R}}
\def\rn{{{\rr}^n}}
\def\zz{{\mathbb Z}}

\def\zz{{\mathbb Z}}

\def\mm{{\mathbb M}}

\def\az{{\alpha}}

\begin{document}
\title{{Optimal asymptotic expansion of entire solutions to Monge-Amp\`{e}re equation with $C^\alpha$ perturbed periodic data} }

\author[1]{ Shuai Qi }
%\quad %and\quad 
\author[,2]{Jiguang Bao\thanks{Corresponding author.	}}
%\author[,a]{Shuai Qi\corref{cor1}}
%the National Natural Science Foundation of China (No.11771023).}
%\author[2]{Dinghuai Wang }
%\ \ \ and \ \ \ % Yannick Sire 
%\quad and\quad 
\affil[1]{School of Sciences, Qingdao University of Technology, Qingdao, 266520,  China}
\affil[2]{School of Mathematical Sciences, Beijing Normal University, 
	%Laboratory of Mathematics and Complex Systems, Ministry of Education, 	
	Beijing, 100875, China}
%\affil[2]{School of Mathematics and Statistics, Anhui Normal University, Wuhu, 241002, China}
%\cortext[cor1]{Corresponding author: Shuai Qi.	}
\date{}
\maketitle

\renewcommand\thefootnote{}
\footnote{{\it E-mail address:} qshuai@pku.edu.cn(S. Qi), jgbao@bnu.edu.cn(J. Bao). }
\setcounter{footnote}{0}
\renewcommand{\thefootnote}{\fnsymbol{footnote}}
\vspace{-1cm}

{\bf Abstract:} We consider the asymptotic behavior at infinity of solution $u$ to Monge-Amp\`{e}re equation $\det(D^2u)=f$ in $\rn$, where 
$f$ is a perturbation of a periodic function and is only assumed to be H\"{o}lder continuous, compared to the previous work that $f$ is at least 
$C^{1,\az}$. The consequence established in this paper, by a nonlocal method, is that the difference between $u$ and a quadratic 
polynomial is asymptotically close to a periodic function.

%{\marginpar{ \cite{tz}摘要中作者陈述其$f\in C^{1,\az}$,但其对$f$的假设要$f\in C^3$}}

%\tableofcontents

%80A25.
\abovedisplayskip=5pt
\abovedisplayshortskip=5pt
\belowdisplayskip=5pt
\belowdisplayshortskip=5pt
\section{Introduction}

Monge-Amp\`{e}re equations are a class of fully nonlinear equations and can be found in several contexts of analysis and geometry. 
The asymptotic behavior of solution to the Monge-Amp\`{e}re equation 
\begin{equation}\label{eqma}
	\det(D^2u)=f,\ \ \ \ \ \text{in}\ \mathbb{R}^n,
\end{equation}
as an extension of Liouville's theorem, has been investigated extensively in the last years, starting with the pioneering works by 
J\"{o}rgens \cite{j}, Calabi \cite{ce} and Pogorelov \cite{p}. They showed that any convex classical solution of \eqref{eqma} with $f\equiv1$ 
must be a quadratic polynomial. 
Caffarelli \cite{ca} extended the result from classical solution to viscosity solution.% in the case $f\equiv1$. 

Caffarelli and Li \cite{cl} 
considered the case $f$ is equal to $1$ outside a compact subset of $\rn$ and Bao, Li and Zhang \cite{blz} focused on a more general situation 
that $f$ is a perturbation of $1$ near infinity. They all proved that any convex viscosity solution $u$ is close to a 
quadratic polynomial at infinity for $n\geq3$ (there is an additional logarithmic term for $n=2$). In \cite{lb}, Bao and Liu concentrated on 
the regularity of $f$ and weakened it from $f\in C^3$ in \cite{blz} to $f\in C^2$ by the spherical harmonic expansion of the solutions of 
nonhomogeneous linearized equations.  Bao and Qi \cite{qb} developed a nonlocal method to obtain the asymptotic behavior of $u$ under 
the hypothesis on regularity that $f$ is just H\"{o}lder continuous. In \cite{jtx}, Jin, Tu and Xiong considered the case that the right-hand side of the equation is a measure.

On the other hand, $f\equiv1$ is also a special periodic function. In the well known work \cite{cl2}, Caffarelli and Li investigated the case 
that $f\in C^\az$ is a positive periodic function and showed that the difference between $u$ and a quadratic polynomial is periodic. Li and Lu considered the 
positive periodic function $f\in L^\infty$ in a following work \cite{ll}. As the counterpart of \cite{blz} in the aperiodic case, Teixeira and Zhang \cite{tz} studied the 
the asymptotic behavior of solution to \eqref{eqma}, in which the right hand side $f\in C^{1,\az}$ is asymptotically close to periodic function. 
They showed that the difference between $u$ and a quadratic polynomial is asymptotically close to a periodic function.

In this paper we discuss the asymptotic behavior of solution $u$ to \eqref{eqma} with a H\"{o}lder continuous term $f$, compared to the condition 
$f\in C^{1,\az}$ in \cite{tz}. More precisely, our assumptions on $f$ are as follows:

Let $f_p$ be a positive periodic function in $\rn$ such that, for some constants $d_0>0$, $\alpha\in(0,1)$ and $a_1,\cdots,a_n>0$ there hold
\begin{equation*}
   \begin{array}{rl}
   	   &d_0^{-1}\leq f_p\leq d_0,\ \ \ [f_p]_{C^\alpha(\rn)}\leq d_0,\\
   	   &f_p(x+a_ie_i)=f_p(x),\ \ \ \forall x\in\rn,
   \end{array}
\end{equation*}
where $e_i=(0,\cdots,0,1,0,\cdots,0)$ denotes the $i$-th standard coordinate vector in $\rn$, $i=1,\cdots,n$.

%We suppose that $f\in C^\alpha(\rn)$ is asymptotically close to $f_p$ in the following sense: 
%\vspace{-0.2cm}
%\begin{itemize}
%	%\item[(H1)] The function $f\in C(\mathbb{R}^n)$ is positive. 
%	\item[(H)] There exist $d_1>0$ and $\beta>2$ such that 
%	\begin{equation*}
%		d_1^{-1}\leq f\leq d_1,\ \ \ \|f\|_{C^\alpha(\rn)}\leq d_1,
%	\end{equation*} 
%	and, for $|x|\geq1$, 
%	\begin{equation}\label{eqfcdinfi}
%		\begin{aligned}
%			|x|^{\be}\left|\left(f^{\frac{1}{n}}-f_p^{\frac{1}{n}}\right)(x)\right|+
%			|x|^{\be+\alpha}\left[\left(f^{\frac{1}{n}}-f_p^{\frac{1}{n}}\right)(x)\right]_{C^\az(\overline{B_{\frac{|x|}{2}}(x)})}\leq d_1.
%		\end{aligned}
%	\end{equation}
%\end{itemize}

We suppose that $f$ is asymptotically close to $f_p$ in the following sense: 
\vspace{-0.2cm}
\begin{itemize}
	%\item[(H1)] The function $f\in C(\mathbb{R}^n)$ is positive. 
	\item[(H)] The function $f\in C(\rn)\cap C_{loc}^{\alpha}({\rn\setminus\overline{\O}})$ for 
	some bounded open subset $\O\subset\rn$ and satisfies, for $x\in\rn\setminus\{0\}$, 
	\begin{equation*}\label{eqfcdinfi}
		\begin{aligned}
			|x|^{\be}\left|\left(f-f_p\right)(x)\right|+
			|x|^{\be+\alpha}\left[f-f_p\right]_{C^\az(\overline{B_{\frac{|x|}{2}}(x)})}\leq d_1
		\end{aligned}
	\end{equation*}
	for some constants $\be>2$ and $d_1>0$.
\end{itemize}
\vspace{-0.2cm}
Let $\mm^{n\times n}$ be the set of the real valued, $n\times n$ matrices and 
\begin{equation*}
	\A:=\left\{A\in\mm^{n\times n}:\ A\text{ is symmetric, positive definite and }\det(A)=
	\int\hspace{-1.05em}-_{\raisebox{-5pt}{\tiny$\prod_{1\leq i\leq n}[0,a_i]$}}fdx\right\}.
\end{equation*}

Our main result is
\begin{theorem}\label{thm1}
	Let $u\in C_{loc}(\rn)$ be a convex viscosity solution of \eqref{eqma}, where $f:\rn\rightarrow\rr$ satisfies 
	{\textnormal{(H)}}. If $n\geq3$, then there exist $b\in\rn$, $A\in\A$ and $u_p\in C^{2,\alpha}(\rn)$ such that
	\begin{equation*}
		\begin{aligned}
		\left|u(x)-\left(\frac{1}{2}x'Ax+b\cdot x+u_p(x)\right)\right|\leq C_0(1+|x|)^{2-\min\{n,\beta\}},\ \ \ \forall\ x\in\rn,
		\end{aligned}
	\end{equation*}
	where $u_p$ is $a_i$-periodic in the $i$-th variable, $i=1,\cdots,n$,   
    and the constant $C_0>0$ depends on $d_0$, $d_1$, $n$, $\be$, $\az$ and $a_1,\cdots,a_n$.
\end{theorem}

\begin{remark}\label{remcf}
	It is clear that $f\in C^\az(\rn\setminus\O)$. Moreover, the viscosity solution $u$ in Theorem \ref{thm1} 
	is really $C^{2,\alpha}$ near infinity, see Section 3.
\end{remark}

%\begin{remark}\label{remr}
%	The viscosity solution $u$ in Theorem \ref{thm1} is really $C^{2,\alpha}$ near infinity and we can obtain this result by repeating the 
%	process in \cite[Lemma 3.2]{qb}.
%\end{remark}

\begin{remark}\label{rembeta}
	The assumption $\be>2$ is essentially optimal, one can see the example in \cite{blz,tz}.
\end{remark}

The major difficulties in this paper are caused by the lower regularity and the aperiodicity of $f$. To get the asymptotic behavior near infinity, 
the authors both in \cite{cl2} and \cite{tz} paid much attention on the second-order increment of $u$. Their vital step is deducing a suitable 
equation of the increment. More specifically, in \cite{cl2}, although $f$ is just H\"{o}lder continuous they obtained a homogeneous equation of the increment 
since $f$ is periodic, while in \cite{tz} the authors got a nonhomogeneous equation in which the nonhomogeneous term can be controlled since the higher 
regularity of $f$. 
In the case $f\in L^\infty$ is periodic \cite{ll}, the authors established the asymptotic results by following the main ideas in 
\cite{cl2} with additional approximation arguments. 
It is not the case in our paper and the second-order increment is not an appropriate object for us. We will propose a novel 
fractional order increment and deduce a nonhomogeneous equation of it. Then we employ the nonlocal ideas in \cite{qb} and refine the 
nonlocal method to measure the decay of the nonhomogeneous term. The asymptotic behavior is finally established with the aid of the fractional order increment.

%In this paper, a non-local method is developed to prove Theorem \ref{thm1}. 
%As mentioned above, the right term $f$ is at least $C^2$ in the 
%previous work and, of course, the derivative of $f$ plays a crucial role in their proof. Our method is essentially different 
%since the $f$ in this paper is just H\"{o}lder continuous and one can not expect the existence of its derivative. We introduce 
%some non-local objects and develop a series of non-local arguments to deal with the problems caused by the insufficient smoothness 
%of $f$.

\subsection*{Notations and the structure of the article}

Throughout the paper, we use the following notations:

The point $x\in\rn$ will also be written as $x=(x_1,x')=(x_1,\cdots,x_n)$. The notation $|x|$ denotes the Euclidean norm of $x$.
$B_r(x)\subset\rn$ denotes the ball centered at $x$ with radius $r$. 
We drop the center if it coincides with the origin, i.e., $B_r=B_r(0)$.

The identity matrix is denoted by $I$. 
For an $n\times n$ matrix $B=(b_{ij})_{n\times n}$, $1\leq i,j\leq n$, $cof_{ij}B$ stands for the algebraic cofactor of 
element $b_{ij}$. We will always assume that $i,j\in\mathbb{N}_+$ and $i,j\leq n$ in this paper 
except for additional explanation.

%The gradient of a function $u$ is denoted by $Du$. $\Delta$ is the Laplacian. 

Given a function $u:\rn\rightarrow\rr$ and a point $x\in\rn$, we denote by $Du(x)$ and $\Delta u(x)$ the gradient  vector and 
the  Laplacian respectively.  
%\begin{equation*}
%	\partial_iu(x)=\frac{\partial u(x)}{\partial x^{(i)}},\ \ \ \ 
%	\partial_{ij}u(x)=\frac{\partial^2u(x)}{\partial x^{(i)}\partial x^{(j)}}\ \ \ \ \text{for }1\leq i,j\leq n
%\end{equation*}
%and define 
%\begin{equation*}
%	D^\gamma u(x):=\partial^{\gamma_1}_{x^{(1)}}\cdots\partial^{\gamma_n}_{x^{(n)}}u(x),\ \ \ \ \ 
%	D^ku(x):=\{D^\gamma u(x):\ |\gamma|=k\}.
%\end{equation*}
%Assigning some ordering to the various partial derivatives, we can regard $D^ku(x)$ as a point in $\mathbb{R}^{n^k}$ 
%and write 
%\begin{equation*}
%	\left|D^ku(x)\right|:=\left(\sum_{|\gamma|=k}\left|D^\gamma u(x)\right|^2\right)^{1/2}. 
%\end{equation*}
%Clearly, $D^0u(x)=u(x)$. For $k=1$ and $k=2$, we regard the elements of $Du(x):=D^1u(x)$ and $D^2u(x)$ as being 
%arranged in the following ways: 
%\begin{gather*}
%	Du(x):=(\partial_1u(x),\cdots,\partial_nu(x))=gradient\ vector,\\
%	D^2u(x):=(\partial_{ij}u(x))_{n\times n}=Hessian\ matrix.
%\end{gather*}
%We denote by $\Delta$ the Laplacian.

Let $k\in\mathbb{N}$, $\kappa\in(0,1)$ and $\Omega\subset\rn$ be open. For multi-indices 
$\gamma=(\gamma_1,\cdots,\gamma_n)\in\mathbb{N}^n$, we let $|\gamma|$ denote the sum of its components. 
The space $C^k(\overline{\Omega})$ consists of functions $u: \overline{\Omega}\rightarrow\rr$, which admit derivatives 
up to order $k$, such that %$u$ and its derivatives up to order $k$ are continuous and bounded on $\overline{\Omega}$. 
\begin{equation*}
	\|u\|_{C^{k}(\overline{\Omega})}:=\sup_{\substack{ \gamma\in\mathbb{N}^n\\|\gamma|\leq k}}
	\sup_{x\in{\Omega}}\left|D^\gamma u(x)\right|<+\infty.%,\ \ \ \ \forall u\in C^{k}(\overline{\Omega}),
\end{equation*}
The H\"{o}lder space $C^{k,\kappa}(\overline{\Omega})$ consists of function $u\in C^k(\overline{\Omega})$ satisfying 
\begin{gather*}
	\|u\|_{C^{k,\kappa}(\overline{\Omega})}:=\|u\|_{C^{k}(\overline{\Omega})}+[u]_{C^{k,\kappa}(\overline{\Omega})}<+\infty, 
\end{gather*} 
where the seminorm 
%The $\beta$-H\"{o}lder seminorm is defined by
\begin{equation*}
	[u]_{C^{k,\kappa}(\overline{\Omega})}:=\sup_{\substack{ \gamma\in\mathbb{N}^n\\|\gamma|=k}}
	\sup_{\substack{x,y\in{\Omega} \\ x\neq y}}\frac{|D^\gamma u(x)-D^\gamma u(y)|}{|x-y|^\kappa}.
\end{equation*}
%We say $u\in C^{\beta}(\overline{\Omega})$ if 
%\begin{equation*}
%	\|u\|_{C^\beta(\overline{\Omega})}:=\|u\|_{C(\overline{\Omega})}+[u]_{C^\beta(\overline{\Omega})}<+\infty.
%\end{equation*}
%We write $C^{0}(\overline{\Omega})$, $C^{0,\beta}(\overline{\Omega})$ as $C(\overline{\Omega})$, 
%$C^{\beta}(\overline{\Omega})$, respectively.
A function $u\in C^{k,\kappa}_{loc}({\Omega})$ if $u\in C^{k,\kappa}(K)$ for all compact 
$K\subset\Omega$. %Clearly, $C^\beta(\overline{\Omega})$ coincides with $C^{0,\beta}(\overline{\Omega})$. 
For convenience, we may write $C^{0}(\overline{\Omega})$, $C^{k,\kappa}(\overline{\Omega})$ as $C(\overline{\Omega})$, 
$C^{k+\kappa}(\overline{\Omega})$,  respectively. 

%If $\det(D^2v(x))=f(x)$ in the classical sense for functions $v$ and $f$, then, by the property of determinant, we have 
%\begin{equation*}
%	\sum_{1\leq i,j\leq n}cof_{ij}(D^2v(x))\partial_{ij}v(x)=nf(x).
%\end{equation*}
%Let $a_{ij}(x):=n^{-1}cof_{ij}(D^2v(x))$, then 
%\begin{equation}\label{notaeq}
%	\sum_{1\leq i,j\leq n}a_{ij}(x)\partial_{ij}v(x)=f(x).
%\end{equation}
%For convenience, we will omit the notation $\sum$ and write \eqref{notaeq} as $a_{ij}(x)\partial_{ij}v(x)=f(x)$ in this article. 

The paper is organized as follows: In Section 2, we introduce the fractional Laplacian and a series of relevant properties,  
which will be applied in the sequel. In Section 3, we deduce the equation for the fractional order increment (Lemma \ref{lemeqv}) and then 
show the crucial estimate (Lemma \ref{lemeF}) by a nonlocal method. Finally, Theorem \ref{thm1} is proved in Section 4.

\section{Fractional Laplacian and its properties}
In this section we introduce the fractional Laplacian and list some crucial properties established in \cite{qb,s}.

%Based on these results, we can use the fractional Laplacian to show Theorem \ref{thm1}.

The fractional Laplacian for a function $u:\rn\rightarrow\rr$ with the parameter $s\in(0,1)$ is defined as
\begin{equation*}
	(-\Delta)^su(x)=c_{n,s} \textnormal{P.V.}\int_{\rn}\frac{u(x)-u(y)}{|x-y|^{n+2s}}dy,
\end{equation*}
where $c_{n,s}$ is the normalization constant and P.V. stands for the Cauchy principal value. 
If $u$ is $C^{2s+\ep}$ at a point $x_0\in\rn$ for some $\ep>0$ and satisfies %{\marginpar{这里$\ep$表示一个小正数,换成$\az$不太合适}}
\begin{equation*}\label{eqfwdc}
	\|u\|_{L_s(\rn)}:=\int_{\rn}\frac{|u(y)|}{1+|y|^{n+2s}}dy<+\infty,
\end{equation*}
then $(-\Delta)^su(x_0)$ can be calculated classically. 

The fundamental solution of $(-\Delta)^s$ (see \cite{b,s}) is 
\begin{equation*}
	\Phi_s(x):=c_{n,-s}\frac{1}{|x|^{n-2s}},\ \ \ \ n>2s.
\end{equation*}
As in \cite{s}, we use the notation 
\begin{equation*}
	(-\Delta)^{-s}u(x)=c_{n,-s}\int_{\rn}\frac{u(y)}{|x-y|^{n-2s}}dy, 
\end{equation*}
which is really the Riesz potential of $u$. If $u\in L^p(\rn)$ and $2sp>n$, then $(-\Delta)^{-s}u\in C(\rn)$, see \cite{ah,dpv}. 
%More details about fractional Laplacian can be found in \cite{clm,s}.

Now we list some important properties for $(-\Delta)^s$ and $(-\Delta)^{-s}$, which will 
be frequently used in our paper. We first summarize \cite[Propositions 2.5-2.7]{s} as a following version.   
\begin{proposition}\label{prorf}
	Let $s\in(0,1)$ and $u\in C^{k,\alpha}(\rn)$ for $\alpha\in(0,1)$ and $k=0,1$.
	\vspace{-0.2cm}
	\begin{itemize}
		\setlength{\itemsep}{-3pt}
		\item If $\alpha>2s$, then $(-\Delta)^su\in C^{k,\alpha-2s}(\rn)$ and 
		\begin{equation*}
			\left[(-\Delta)^su\right]_{C^{k,\alpha-2s}(\rn)}\leq C[u]_{C^{k,\alpha}(\rn)}.
		\end{equation*}
		\item If $\alpha<2s<1+\alpha$ and $k=1$, then $(-\Delta)^su\in C^{\alpha-2s+1}(\rn)$ and 
		\begin{equation*}
			\left[(-\Delta)^su\right]_{C^{\alpha-2s+1}(\rn)}\leq C[u]_{C^{1,\alpha}(\rn)}.
		\end{equation*}
	\end{itemize}
	\vspace{-0.2cm}
	The constant $C$ depends only on $\alpha,\ s$ and $n$.
\end{proposition}

Next we introduce two properties established in \cite{qb}. For reading simplicity, the details of the proof are provided here. 
The first is a decay result for fractional Laplacian (\cite[Lemma 2.5]{qb}).
\begin{proposition}\label{prodf}
	Let $\alpha\in(0,1)$ and $s\in(0,\frac{\alpha}{2})$. If $u\in C^{\alpha}(\rn)$ satisfies
	\begin{equation}\label{eqcfd2.5}
		|u(x)|\leq c'|x|^{-\kappa},\ \ [u]_{C^{\alpha}(\overline{{B_{\frac{|x|}{2}}(x)}})}\leq c'|x|^{-\kappa-\alpha},\ \ \ \text{for}\ |x|>R'
	\end{equation}
	for some positive constants $c',\kappa$ and $R'$, then there exists $C_1>0$ such that 
	\begin{equation*}
		|(-\Delta)^su(x)|\leq 
		\left\{
		\begin{aligned}
			C_1|x|^{-\min\{\kappa,n\}-2s},\ \ \ \kappa\neq n,\\
			C_1|x|^{-n-2s}(\ln|x|),\ \ \ \kappa=n,
		\end{aligned}
		\right.\ \ \ \ \ \ \text{for}\ |x|>2R',
	\end{equation*}    
	where the constant $C_1$ depends only on $n,\ s,\ \kappa,\ c'$, $R'$ and $\|u\|_{C(\rn)}$.
\end{proposition}

\noindent{\it Proof.} Without loss of generality, we assume $R'>1$. Define for $|x|>2R'$ that 
\begin{align*}
	&A_1:=\{y\in\rn:\ |y|\leq\frac{|x|}{2}\},\\
	&A_2:=\{y\in\rn:\ |x-y|\leq\frac{|x|}{2}\},\\
	&A_3:=\rn\setminus\left(A_1\cup A_2\right).
\end{align*}
Then 
\begin{equation*}
	(-\Delta)^su(x)%=c_{n,s}\int_{\rn}\frac{u(x)-u(y)}{|x-y|^{n+2s}}dy
	=c_{n,s}\int_{A_1\cup A_2\cup A_3}\frac{u(x)-u(y)}{|x-y|^{n+2s}}dy.
\end{equation*}

Note that $|x-y|\geq\frac{|x|}{2}\geq|y|$ in $A_1$ and $\frac{|x|}{2}>1$, we write  
\begin{align*}
	\int_{A_1}\frac{|u(x)-u(y)|}{|x-y|^{n+2s}}dy
	&\leq\frac{C}{|x|^{n+2s}}\int_{B_{\frac{|x|}{2}}}|u(x)-u(y)|dy\\
	&\leq\frac{C}{|x|^{n+2s}}\left(\int_{B_{R'}}+\int_{B_{\frac{|x|}{2}}\setminus B_{R'}}\right)|u(x)-u(y)|dy.
\end{align*}
Since $u$ is bounded, one can find a constant $C'>0$, which may depend on $\|u\|_{C(\rn)}$ and $R'$ such that   
\begin{equation*}
	\int_{B_{R'}}|u(x)-u(y)|dy%\leq\int_{B_1}|u(x)|dy+\int_{B_1}|u(y)|dy\leq C_n(R')^{-\kappa}+C
	\leq C'.
\end{equation*}
On the other hand, equation \eqref{eqcfd2.5} yields that 
\begin{align*}
	\int_{B_{\frac{|x|}{2}}\setminus B_{R'}}|u(x)-u(y)|dy\leq C\int_{B_{\frac{|x|}{2}}\setminus B_{R'}}|y|^{-\kappa}dy\leq
	\left\{
	\begin{array}{ll}
		C|x|^{n-\kappa},&\ \ \ \kappa\neq n,\\
		C\ln|x|,&\ \ \ \kappa=n.
	\end{array}
	\right.
\end{align*}
%{\color{red} I think it should be $C|x|^{n-\kappa}+C$ in the case $\kappa\neq n$ but not $C|x|^{n-\kappa}$. 
%	The $+C$ is the value at $|y|=1$. If we just write $C|x|^{n-\kappa}$ in this case, 
%	then the constant $C$ here may depend on $x$ since $n-\kappa$ may be negative and $|x|$ is large.}\\
Thus 
\begin{align*}
	\int_{B_{\frac{|x|}{2}}}|u(x)-u(y)|dy
	&\leq\left\{
	\begin{array}{ll}
		C(|x|^{n-\kappa}+1),&\ \ \ \kappa\neq n,\\
		C(\ln|x|+1),&\ \ \ \kappa=n.
	\end{array}
	\right.
\end{align*}
Note that 
\begin{align*}
	&(|x|^{n-\kappa}+1)\cdot{|x|^{-n-2s}}\leq|x|^{-\min\{\kappa,n\}-2s},\\
	&(\ln|x|+1)\cdot{|x|^{-n-2s}}\leq|x|^{-n-2s}(\ln|x|),
\end{align*}
we obtain    
\begin{align*}
	\int_{A_1}\frac{|u(x)-u(y)|}{|x-y|^{n+2s}}dy
	&\leq\left\{
	\begin{aligned}
		C|x|^{-\min\{\kappa,n\}-2s},\ \ \ \kappa\neq n,\\
		C|x|^{-n-2s}(\ln|x|),\ \ \ \kappa=n,
	\end{aligned}
	\right.\ \ \ \text{for}\ |x|>2R'.
\end{align*}
%The estimates in $A_2$ and $A_3$ are the same as in Lemma \ref{lemd1}. The lemma is proved.

In $A_2$ we have $|x-y|\leq\frac{|x|}{2}\leq|y|$. Then \eqref{eqcfd2.5} allows us to calculate that  
\begin{align*}
	\int_{A_2}\frac{|u(x)-u(y)|}{|x-y|^{n+2s}}dy
	\leq\frac{c'}{|x|^{\kappa+\alpha}}\int_{|x-y|\leq\frac{|x|}{2}}\frac{1}{|x-y|^{n+2s-\alpha}}dy
	\leq C|x|^{-\kappa-2s}.
\end{align*}

For $A_3$, we divide it into two parts
\begin{equation*}
	A_3^+:=\{y\in A_3:\ |x-y|\geq|y|\}\ \ \ \text{and}\ \ \ A_3^-:=A_3\setminus A_3^+.
\end{equation*}
Then we have $\frac{|x|}{2}\leq|y|\leq|x-y|$ in $A_3^+$ and hence $|u(x)|+|u(y)|\leq C|x|^{-\kappa}$, which gives 
\begin{align*}
	\int_{A_3^+}\frac{|u(x)-u(y)|}{|x-y|^{n+2s}}dy\leq
	\frac{C}{|x|^{\kappa}}\int_{|y|\geq\frac{|x|}{2}}\frac{1}{|y|^{n+2s}}dy
	\leq C|x|^{-\kappa-2s}.
\end{align*}
Finally, in $A_3^-$ there holds $\frac{|x|}{2}\leq|x-y|\leq|y|$ and thus we have 
\begin{align*}
	\int_{A_3^-}\frac{|u(x)-u(y)|}{|x-y|^{n+2s}}dy\leq
	\frac{C}{|x|^{\kappa}}\int_{|x-y|\geq\frac{|x|}{2}}\frac{1}{|x-y|^{n+2s}}dy
	\leq C|x|^{-\kappa-2s}.
\end{align*}
The lemma follows from the above estimates.
\qed

For the Riesz potential $(-\Delta)^{-s}$, which can be seen as the inverse of $(-\De)^s$, we also has a similar decay ( \cite[Corollary 2.9]{qb}).
\begin{proposition}\label{prodif}
	Let $s\in(0,1)$ and $u\in C(\rn)$ be such that 
	\begin{equation*}
		|u(x)|\leq c'|x|^{-\kappa},\ \ \ \text{for}\ |x|>R'
	\end{equation*}
	for some positive constants $c',\kappa>2s$ and $R'$. Then 
	\begin{equation*}
		|(-\Delta)^{-s}u(x)|\leq %C|x|^{-\kappa+2s},\ \ \ |x|>R.
		\left\{
		\begin{aligned}
			C|x|^{2s-\min\{\kappa,n\}},\ \ \ \kappa\neq n,\\
			C|x|^{2s-n}(\ln|x|),\ \ \ \kappa=n,
		\end{aligned}
		\right.\ \ \ \ \ \ \text{for}\ |x|>2R',
	\end{equation*} 
	where the constant $C>0$ depends only on $n,\ s,\ \kappa,\ c'$, $R'$ and $\|u\|_{C(\rn)}$.   
\end{proposition}

\noindent{\it Proof.} It is easy to check that $u\in L^p(\rn)$ for $p$ large enough, 
thus $(-\Delta)^{-s}u\in C(\rn)$, see \cite[Theorem 1.2.4]{ah}. 
Let $A_1,\ A_2$ and $A_3$ be as in Proposition \ref{prodf}, then 
\begin{equation*}
	(-\Delta)^{-s}u(x)=c_{n,-s}\int_{\rn}\frac{u(y)}{|x-y|^{n-2s}}dy
	=c_{n,-s}\int_{A_1\cup A_2\cup A_3}\frac{u(y)}{|x-y|^{n-2s}}dy.
\end{equation*}
The estimate in $A_1$ is similar to that of Proposition \ref{prodf}, one just needs to replace $s$ by $-s$. 
In $A_2$ we have $|y|\geq\frac{|x|}{2}\geq|x-y|$ and thus 
\begin{align*}
	\int_{A_2}\frac{|u(y)|}{|x-y|^{n-2s}}dy
	&\leq\frac{C}{|x|^{\kappa}}\int_{B_{\frac{|x|}{2}}(x)}\frac{1}{|x-y|^{n-2s}}dy
	%&\leq\frac{C}{|x|^{\kappa}}\int_{B_{\frac{|x|}{2}}(x)}\frac{1}{|x-y|^{n-2s}}dy\\
	\leq C|x|^{2s}\cdot|x|^{-\kappa}\leq C|x|^{2s-\kappa}.
\end{align*}
To compute the integral in $A_3$, we also use the notations $A_3^+$ and $A_3^-$ as in Proposition \ref{prodf}. Note that we have   
$|u(y)|\leq c'|y|^{-\kappa}$ and $|x-y|\geq|y|\geq\frac{|x|}{2}$ in $A_3^+$, which yield  
\begin{align*}
	\int_{A_3^+}\frac{|u(y)|}{|x-y|^{n-2s}}dy\leq
	C\int_{|y|\geq\frac{|x|}{2}}\frac{1}{|y|^{n-2s+\kappa}}dy
	\leq C|x|^{2s-\kappa}
\end{align*}
since $\kappa>2s$. In $A_3^-$ the inequalities $\frac{|x|}{2}\leq|x-y|\leq|y|$ hold, thus 
$|u(y)|\leq c'|y|^{-\kappa}\leq c'|x-y|^{-\kappa}$ and 
\begin{align*}
	\int_{A_3^-}\frac{|u(y)|}{|x-y|^{n+2s}}dy\leq
	C\int_{|x-y|\geq\frac{|x|}{2}}\frac{1}{|x-y|^{n-2s+\kappa}}dy
	\leq C|x|^{2s-\kappa}.
\end{align*}
The proof is completed. 
\qed

\section{The study of the fractional order increment}
The aim of this section is to investigate the fractional order increment. For convenience we assume $a_1=\cdots=a_n=1$ and 
$\int\hspace{-0.9em}-_{\raisebox{-0.2pt}{\tiny$\prod_{1\leq i\leq n}[0,a_i]$}}fdx=1$ since the 
Monge-Amp\`{e}re equation is invariant under affine transformation.

We summarize from the first two steps in \cite{tz} the following two propositions.

\begin{proposition}\label{pro1}
	Under the assumptions in Theorem \ref{thm1}, there exists a linear transform $T$ satisfying $\det T=1$ such that
	\begin{equation*}
		\left|\left(u\circ T\right)(x)-\frac{1}{2}|x|^2\right|\leq c_0|x|^{2-\theta},\ \ \ \text{for}\ |x|\geq1
	\end{equation*}
    for some $c_0>0$ and $\theta>0$, depending only on $n,\ \beta,\ d_0$ and $d_1$.
\end{proposition}

This is a important result in studying the asymptotic behavior and one can find the counterparts in \cite{blz,cl2,lb,qb,tz}.
In fact, such a result allows us to get the regularity of $u$. We always write $v=u\circ T$ from now. According to Proposition \ref{pro1}, 
there exist $c_0>0$ and $\theta>0$ such that
\begin{equation*}\label{eqv1}
	\left|v(x)-\frac{1}{2}|x|^2\right|\leq c_0|x|^{2-\theta},\ \ \ \ \text{for}\ |x|\geq1.
\end{equation*}

The equation for $v$ is 
\begin{equation}\label{eqv}
	\det(D^2v(x))=f_T(x),\ \ \ \ \forall x\in\rn,
\end{equation}
where $f_T(x)=f(T(x))$; correspondingly, we let $(f_p)_T(x)=f_p(T(x))$.
The regularity of $v$ is as follows:
\begin{proposition}\label{pro2}
	The function $v\in C^2$ and there exist positive $c_1$ and $c_2$, depending only on $n,\ d_0,\ d_1,\ \be,$ and $\az$ such that 
	\begin{equation*}
		c_1I\leq D^2v(x)\leq c_2I,\ \ \ \ \forall x\in\rn.
	\end{equation*}  
\end{proposition}

From Remark \ref{remcf} we see that $f\in C^\az(\rn\setminus\O)$. By using the classical results in \cite{ca1,c,fjm} and 
repeating the process in \cite[Lemma 3.2]{qb}, we conclude that the viscosity solution $u$ in Theorem \ref{thm1} 
is really $C^{2,\alpha}$ near infinity, i.e., $u\in C^{2,\az}_{loc}(\rn\setminus\overline{B_{2R}})$ for some $R$ large. 
Thus we can use \cite[Theorem 3.4]{y} to redefine the value of $u$ inside $B_{3R}$ such that the new convex function $\bar{u}$ is locally $C^{2,\az}$ on 
$\rn$. Let $\bar{f}:=\det(D^2\bar{u})$, then $\bar{f}\in C^\alpha(\rn)$ coincides with $f$ outside $B_{3R}$. 
Therefore, from now, we assume that $u\in C^{2,\az}_{loc}(\rn)$ and $f\in C^\az(\rn)$.
%This means that we can prove Theorem \ref{thm1} under the additional hypotheses $u\in C^{2,\az}_{loc}(\rn)$ and $f\in C^\az(\rn)$. 

\begin{remark}\label{rem1}
	By assumption (H) and the hypothesis $f\in C^\az(\rn)$ we really have 
	\begin{equation*}
		\begin{array}{rl}
			\left|\left(f-f_p\right)(x)\right|&\leq C(1+|x|)^{-\be},\\
			\left[\left(f-f_p\right)\right]_{C^\az(\overline{B_{\frac{|x|}{2}}(x)})}&\leq C(1+|x|)^{-\be-\az},
		\end{array}
		\ \ \ \ \forall |x|\geq1,
	\end{equation*}
	since $|x|\leq1+|x|\leq2|x|$ when $|x|\geq1$. Note that $f,f_p\in C^\az(\rn)$, we also have 
	\begin{equation*}
		\begin{array}{rl}
			\left|\left(f-f_p\right)(x)\right|&\leq C(1+|x|)^{-\be},\\
			\left[\left(f-f_p\right)\right]_{C^\az(\overline{B_{1}(x)})}&\leq C(1+|x|)^{-\be-\az},
		\end{array}
		\ \ \ \ \forall |x|\leq1.
	\end{equation*}
	The constant $C$ depends only on $d_0$, $d_1$, $n$, $\be$ and $\az$. %$\|f_p\|_{C^\az(\rn)}$ and $\|f\|_{C^\az(\rn)}$.
\end{remark}

Now we define an important set 
\begin{equation*}
	\begin{aligned}
		E:=\big\{k_1p_1+\cdots+k_np_n:\ k_1,\cdots,k_n\in\zz\},
		          %&\ \ \ \ \ (f_p)_T(x+p_i)=(f_p)_T(x),\ \forall x\in\rn,\ 1\leq i\leq n\big\}.
	\end{aligned}
\end{equation*}
where $p_i:=T^{-1}e_i\in\rn$ satisfying
\begin{equation*}
	(f_p)_T(x+p_i)=(f_p)_T(x),\ \ \ \ \forall x\in\rn,\ 1\leq i\leq n.
\end{equation*}

For $e\in\rn\setminus\{0\}$ and $s\in(0,1)$, we introduce the following $2s$-order increment for $v$:
\begin{equation*}
	\De_e^sv(x):=\frac{v(x+e)+v(x-e)-2v(x)}{|e|^{2s}}
\end{equation*}
We concentrate on the $2s$-order increment from now. We first derive the inequality for $\De_e^sv$ from the concavity of $\det^{\frac{1}{n}}$. 
\begin{lemma}\label{lemeqv}
	Let $E$, $v$ and $f_T$ be as above. Then, for $s\in(0,\frac{\az}{2})$, 
	\begin{equation*}
		a_{ij}(x)\pat_{ij}(\De_e^sv(x))\geq F_s(x),\ \ \ \ \forall x\in\rn,\ e\in E\setminus\{0\},
	\end{equation*}
	where $a_{ij}(x)=cof_{ij}(D^2v(x))$ and 
	\begin{equation*}
		F_s(x):=n\left(\det(D^2v(x))^{\frac{n-1}{n}}\right)\left(\frac{f_T^{\frac{1}{n}}(x+e)+f_T^{\frac{1}{n}}(x-e)-2f_T^{\frac{1}{n}}(x)}{|e|^{2s}}\right).
	\end{equation*}
\end{lemma}

\noindent{\it{Proof. }} For convenience we let $G=\det^{\frac{1}{n}}$ and 
\begin{equation*}
	G_{ij}(B)=\frac{\pat_{ij}G}{\pat b_{ij}}(B),\ \ \ \ \text{for}\ B=(b_{ij})_{n\times n}.
\end{equation*}
Then by direct calculation one obtains 
\begin{equation}\label{eqlemeqv1}
	G_{ij}(D^2v(x))=\frac{1}{n}\left(\det(D^2v(x))\right)^{\frac{1-n}{n}}cof_{ij}(D^2v(x)).
\end{equation}
Define for $e\in E$ that 
\begin{equation*}
	w(x)=\frac{v(x+e)+v(x-e)}{2}.
\end{equation*}
Then the concavity of $G$ gives 
\begin{equation}\label{eqlemeqv2}
	\begin{aligned}
		G(D^2w(x))&\geq\frac{1}{2}\left(G(D^2v(x+e))+G(D^2v(x-e))\right)\\
		&=\frac{1}{2}\left(f_T^{\frac{1}{n}}(x+e)+f_T^{\frac{1}{n}}(x-e)\right).
	\end{aligned}
\end{equation}
On the other hand, we can get from the concavity of $G$ that 
\begin{equation}\label{eqlemeqv3}
	\begin{aligned}
		G(D^2w(x))&\leq G(D^2v(x))+G_{ij}(D^2v(x))\pat_{ij}(w(x)-v(x))\\
		&=f_T^{\frac{1}{n}}(x)+G_{ij}(D^2v(x))\pat_{ij}(w(x)-v(x)).
	\end{aligned}
\end{equation}
Note that $|e|^{2s}\De_e^sv(x)=2\left(w(x)-v(x)\right)$, the lemma follows from \eqref{eqlemeqv1}-\eqref{eqlemeqv3}.
\qed

Note the inequality for $\De_e^sv$ is nonhomogeneous, we need to establish an estimate to measure its decay near infinity. It is 
crucial in studying the asymptotic behavior. Let us observe the decay of $F_s$. 

Since $e\in E\setminus\{0\}$, we have $(f_p)_T^{\frac{1}{n}}(x+e)+(f_p)_T^{\frac{1}{n}}(x-e)-2(f_p)_T^{\frac{1}{n}}(x)=0$. Thus we rewrite $F_s$ as 
\begin{equation*}
	\begin{aligned}
		F_s(x)=n\left(\left(f_T(x)\right)^{\frac{n-1}{n}}\right)\left(\frac{g(x+e)-g(x)}{|e|^{2s}}-\frac{g(x)-g(x-e)}{|e|^{2s}}\right),
	\end{aligned}
\end{equation*}
where $g(x):=f_T^{\frac{1}{n}}(x)-(f_p)_T^{\frac{1}{n}}(x)$.
One can apply assumption (H) and Remark \ref{rem1} to get that $g\in C^\az(\rn)$ and 
\begin{equation*}
	\begin{aligned}
		|g(x)|&\leq C(1+|x|)^{-\beta},\ \ \ \ \ \ \forall x\in\rn,\\
		[g]_{C^\az(\overline{B_{1}(x)})}&\leq C(1+|x|)^{-\be-\az},\ \ \ \forall |x|\leq1,\\
		[g]_{C^\az(\overline{B_{\frac{|x|}{2}}(x)})}&\leq C(1+|x|)^{-\be-\az},\ \ \ \forall |x|\geq1.
	\end{aligned}
\end{equation*}
We consider the following three cases of $x\in\rn$ and $e\in E\setminus\{0\}$:
\begin{equation*}
	\begin{array}{ll}
		\textnormal{Case I:}&\ |x|\leq2,\\
		\textnormal{Case II:}&\ |x|>2,\ |e|\leq\frac{|x|}{2},\\
		\textnormal{Case III:}&\ |x|>2,\ |e|>\frac{|x|}{2},\\
	\end{array}
\end{equation*}

In \textnormal{Case II}, obviously, one can check that $x+e,x-e\in\overline{B_{\frac{|x|}{2}}(x)}$ and hence 
\begin{align*}
	&\frac{|g(x+e)-g(x)|}{|e|^{2s}}+\frac{|g(x)-g(x-e)|}{|e|^{2s}}\\
	\leq&2[g]_{C^\az(\overline{B_{\frac{|x|}{2}}(x)})}|e|^{\az-2s}
	\leq C(1+|x|)^{-\be-\az}(1+|x|)^{\az-2s}=C(1+|x|)^{-\be-2s}.
\end{align*}
Therefore, 
\begin{equation}\label{eqfff}
	|F_s(x)|\leq C(1+|x|)^{-\be-2s}
\end{equation}

In \textnormal{Case I}, we have
\begin{equation*}
	\frac{|g(x+e)-g(x)|}{|e|^{2s}}+\frac{|g(x)-g(x-e)|}{|e|^{2s}}\leq2[g]_{C^\az(\rn)}|e|^{\az-2s}\leq C
\end{equation*}
for $|e|\leq1$ since $g$ is $C^\az$ in $\rn$. If $|e|>1$, then the boundedness of $g$ implies that 
\begin{equation*}
	\frac{|g(x+e)-g(x)|}{|e|^{2s}}+\frac{|g(x)-g(x-e)|}{|e|^{2s}}\leq4\|g\|_{C(\rn)}\leq C.
\end{equation*}
Thus \eqref{eqfff} holds.
%{\marginpar{$|e|$有正下界,但是Case I这里这样写证明可以让常数$C$不依赖于$|e|$的下界.}}
\begin{remark}\label{remE}
	Note that the set $E$ consists of all the period of function $(f_p)_T$, we can find a constant $c>0$ such that 
	\begin{equation*}
		\inf_{e\in E\setminus\{0\}}|e|\geq c>0,
	\end{equation*}
	since the matrix $T$ in Proposition \ref{pro1} satisfies $\det T=1$.
    Thus the proof in \textnormal{Case I} can be simplified by using this property.
    In this case the constant $C$ may depend on $E$ but it has no influence on Theorem \ref{thm1}.
\end{remark}

For \textnormal{Case III}, we consider the following two parts: 
\begin{gather*}
	\Omega _1:=\left\{|x|>2,\ |e|>\frac{|x|}{2}:\ |x-e|\geq\frac{|x|}{4}\ \text{and }|x+e|\geq\frac{|x|}{4}\right\},\\
	\Omega_2:=\left\{|x|>2,\ |e|>\frac{|x|}{2}:\ |x-e|<\frac{|x|}{4}\ \text{or }|x+e|<\frac{|x|}{4}\right\}.
\end{gather*}
In $\Omega_1$, we have
\begin{gather*}
	|g(x+e)|\leq C(1+|x+e|)^{-\be}\leq C(1+\frac{|x|}{4})^{-\be}\leq C(1+|x|)^{-\be},\\
	|g(x-e)|\leq C(1+|x-e|)^{-\be}\leq C(1+\frac{|x|}{4})^{-\be}\leq C(1+|x|)^{-\be},
\end{gather*}
which yields  
\begin{gather*}
	|g(x+e)-g(x)|\leq C(1+|x|)^{-\be}\ \ \ \text{and}\ \ \ |g(x-e)-g(x)|\leq C(1+|x|)^{-\be}.
\end{gather*}
Since $|x|>2$ and $|e|>\frac{|x|}{2}$, \eqref{eqfff} follows from the fact that 
\begin{equation*}
	|e|^{-2s}\leq C|x|^{-2s}\leq C(1+|x|)^{-2s}.
\end{equation*}

Establishing the decay in $\Omega_2$ is really difficult and we employ a nonlocal method to do it. Finally, we conclude the 
following decay result.

\begin{lemma}\label{lemeF}
	Let $F_s$ be the one obtained in Lemma \ref{lemeqv}. Then we have 
	\begin{equation}\label{eqlemeF1}
		|F_s(x)|\leq 
		\left\{
		\begin{array}{ll}
			C_0(1+|x|)^{-\min\{\be+2s,n\}},&\ \ \ \be+2s\neq n,\ \be\neq n,\\
			C_0(1+|x|)^{-\beta-2s}(\ln(2+|x|)),&\ \ \ \be+2s=n,\\
			C_0(1+|x|)^{-\beta-2s}(\ln(2+|x|)),&\ \ \ \be=n,
		\end{array}
		\right.\ \ \ \ \ 
	\end{equation}
	for all $e\in E\setminus\{0\},\ x\in\rn.$ The constant $C_0>0$ depends only $d_0,\ d_1,\ \be,\ n$ and $s$.
\end{lemma}

To show this lemma, we first establish an auxiliary result.

\begin{lemma}\label{lemeFl}
	Let $s\in(0,\frac{\az}{2})$, $x,e\in\Omega_2$. Define
	\begin{equation}\label{eqpsi}
		\Psi:=
		\left(\frac{1}{|x+e-y|^{n-2s}}-\frac{1}{|x-y|^{n-2s}}\right)-\left(\frac{1}{|x-y|^{n-2s}}-\frac{1}{|x-e-y|^{n-2s}}\right)
	\end{equation}
	for $|y|\geq\frac{|x|}{2}$ and $y\neq x$. Then there exists $C=C(n,s)>0$ such that
	\begin{equation}\label{eqlemeFJ}
		|\Psi|\leq\frac{C|e|^2}{|x-y|^{n+2-2s}}.
	\end{equation}
\end{lemma}

\noindent{\it{Proof. }} Without loss of generality, we assume that $e=(|e|,0,\cdots,0)$. A direct computation shows that 
\begin{equation*}
	\begin{aligned}
		\partial_1|x|^{2s-n}&:=\frac{\partial |x|^{2s-n}}{\partial x_1}=\frac{(2s-n)x_1}{|x|^{n+2-2s}},\\
		\partial_{11}|x|^{2s-n}&:=\frac{\partial^2 |x|^{2s-n}}{\partial (x_1)^2}=(2s-n)\frac{|x|^2+(2s-n-2)x_1^2}{|x|^{n+4-2s}},
	\end{aligned}
\end{equation*}
which allows us to rewrite $\Psi$ as 
%\begin{equation*}
%	\begin{aligned}
	%		\Gamma=&\frac{c_{n,s}}{|e|^{2s}}\int_{\rn}(-\De)^sg(y)\bigg[\int_0^{|e|}\frac{(2s-n)(x_1-y_1+t)}{|(x_1-y_1+t)^2+|x'-y'|^2|^{\frac{n+2-2s}{2}}}dt\\
	%		&\ \ \ \ -\int_0^{|e|}\frac{(2s-n)(x_1-y_1-|e|+t)}{|(x_1-y_1-|e|+t)^2+|x'-y'|^2|^{\frac{n+2-2s}{2}}}dt\bigg]dy\\
	%		=&\frac{c_{n,s}}{|e|^{2s}}\int_{\rn}(-\De)^sg(y)\int_0^{|e|}\int_0^{|e|}\bigg[
	%		\frac{(2s-n)[(x_1-y_1-|e|+t+\tau)^2+|x'-y'|^2]}{|(x_1-y_1-|e|+t+\tau)^2+|x'-y'|^2|^{\frac{n+4-2s}{2}}}\\
	%		&\ \ \ \ +\frac{(2s-n)(2s-n-2)(x_1-y_1-|e|+t+\tau)^2}{|(x_1-y_1-|e|+t+\tau)^2+|x'-y'|^2|^{\frac{n+4-2s}{2}}}\bigg]d\tau dtdy\\
	%		=&\frac{c_{n,s}}{|e|^{2s-2}}\int_{\rn}(-\De)^sg(y)\int_0^{1}\int_0^{1}\bigg[
	%		\frac{(2s-n)[(x_1-y_1+|e|(-1+t+\tau))^2+|x'-y'|^2]}{|(x_1-y_1+|e|(-1+t+\tau))^2+|x'-y'|^2|^{\frac{n+4-2s}{2}}}\\
	%		&\ \ \ \ +\frac{(2s-n)(2s-n-2)(x_1-y_1+|e|(-1+t+\tau))^2}{|(x_1-y_1+|e|(-1+t+\tau))^2+|x'-y'|^2|^{\frac{n+4-2s}{2}}}\bigg]d\tau dtdy\\
	%		\leq&C|e|^{2-2s}\int_{\rn}J\cdot|(-\De)^sg(y)|dy
	%	\end{aligned}
%\end{equation*}
\begin{align*}
	%\begin{aligned}
	\Psi
	%		=&\int_0^{|e|}\frac{(2s-n)(x_1-y_1+t)}{|(x_1-y_1+t)^2+|x'-y'|^2|^{\frac{n+2-2s}{2}}}dt\\
	%		&\ \ \ \ -\int_0^{|e|}\frac{(2s-n)(x_1-y_1-|e|+t)}{|(x_1-y_1-|e|+t)^2+|x'-y'|^2|^{\frac{n+2-2s}{2}}}dt\\
	%		=&\int_0^{|e|}\int_0^{|e|}\bigg[\frac{(2s-n)[(x_1-y_1-|e|+t+\tau)^2+|x'-y'|^2]}{|(x_1-y_1-|e|+t+\tau)^2+|x'-y'|^2|^{\frac{n+4-2s}{2}}}\\
	%		&\ \ \ \ +\frac{(2s-n)(2s-n-2)(x_1-y_1-|e|+t+\tau)^2}{|(x_1-y_1-|e|+t+\tau)^2+|x'-y'|^2|^{\frac{n+4-2s}{2}}}\bigg]d\tau dt\\
	=&|e|^{2}\int_0^{1}\int_0^{1}\bigg[\frac{(2s-n)[(x_1-y_1+|e|(-1+t+\tau))^2+|x'-y'|^2]}{|(x_1-y_1+|e|(-1+t+\tau))^2+|x'-y'|^2|^{\frac{n+4-2s}{2}}}\\
	&\ \ \ \ +\frac{(2s-n)(2s-n-2)(x_1-y_1+|e|(-1+t+\tau))^2}{|(x_1-y_1+|e|(-1+t+\tau))^2+|x'-y'|^2|^{\frac{n+4-2s}{2}}}\bigg]d\tau dt.
	%\end{aligned}
\end{align*}
Clearly, 
\begin{equation*}
	|\Psi|\leq C|e|^2\int_0^{1}\int_0^{1}\frac{1}{|(x_1-y_1+|e|(-1+t+\tau))^2+|x'-y'|^2|^{\frac{n+2-2s}{2}}}d\tau dt
\end{equation*}
If $|x'-y'|\geq\frac{1}{2}|x-y|$, i.e., $|x_1-y_1|\leq\frac{\sqrt{3}}{2}|x-y|$, then obviously
\begin{equation*}
	|\Psi|\leq C|e|^2\int_0^{1}\int_0^{1}\frac{1}{|x'-y'|^{n+2-2s}}d\tau dt\leq\frac{C|e|^2}{|x-y|^{n+2-2s}}.
\end{equation*}
If $|x_1-y_1|>\frac{1}{2}|x-y|$, then we have 
\begin{equation*}
	\begin{aligned}
		|\Psi|\leq&C|e|^2\int_0^{1}\int_0^{1}\frac{1}{|x_1-y_1+|e|(-1+t+\tau)|^{n+2-2s}}d\tau dt\\
		\leq&\frac{C|e|^2}{|e|^{n+2-2s}}\int_0^{1}\int_0^{1}\frac{1}{|L-1+t+\tau|^{n+2-2s}}d\tau dt,
	\end{aligned}
\end{equation*}
where $L:=\frac{x_1-y_1}{|e|}$. One has $|L-1+t+\tau|\geq\frac{|L|}{2}$ when $|L|\geq8$. Thus, \eqref{eqlemeFJ} follows from  
\begin{equation*}
	|\Psi|\leq\frac{C|e|^2}{|e|^{n+2-2s}}\frac{1}{|L|^{n+2-2s}}=\frac{C|e|^2}{|x_1-y_1|^{n+2-2s}}\leq\frac{C|e|^2}{|x-y|^{n+2-2s}}.
\end{equation*}
If $|L|<8$, then $|e|>\frac{|x_1-y_1|}{8}>\frac{|x-y|}{16}$. Thus we have 
\begin{gather*}
	|x\pm e-y|\geq|y|-|x\pm e|\geq\frac{|x|}{4}\geq\frac{|e|}{5}%\ \ \ \text{and}\ \ \ |x\pm e-y|\leq|e|+|x-y|\leq17|e|
\end{gather*}
for $x,e\in\Omega_2$ and $|y|\geq\frac{|x|}{2}$ since $|e|\leq\frac{5}{4}|x|$ here. It follows 
\begin{equation*}
	\frac{1}{|x\pm e-y|^{n-2s}}\leq\frac{C}{|e|^{n-2s}}\leq\frac{C}{|x-y|^{n-2s}}=\frac{C|x-y|^2}{|x-y|^{n+2-2s}}\leq
	\frac{C|e|^2}{|x-y|^{n+2-2s}},
\end{equation*}
which allows us to obtain \eqref{eqlemeFJ} from \eqref{eqpsi}. 
\qed

Next we give a remark on Rolle's theorem, which is important in proving Lemma \ref{lemeF}.

\begin{remark}\label{remrolle}
	Let us recall the Rolle's theorem: Let $z\in C[a,b]$ be differentiable in $(a,b)$ and $z(a)=z(b)$, then there exists a $\xi\in(a,b)$ such that 
	$z'(\xi)=0$. 
	
	Samelson \cite{sa} provided a new proof and concluded that the $\xi$ is not necessarily the absolute (or even a local) maximum or 
	minimum. To enhance readability, we restate the proof in \cite{sa}.
	
	\noindent{\bf{Proof of Rolle's theorem:}} {\rm Define $G(x):=z(x)-z\left(x+\frac{b-a}{2}\right)$, then one has 
	\begin{equation*}
		G(a)=z(a)-z\left(\frac{a+b}{2}\right)\ \ \ \text{and}\ \ \ G\left(\frac{a+b}{2}\right)=z\left(\frac{a+b}{2}\right)-z(b)=-G(a),
	\end{equation*}
	since $z(a)=z(b)$. Thus there exists a $a_1\in[a,\frac{a+b}{2}]$ such that $G(a_1)=0$. Let $b_1=a_1+\frac{b-a}{2}$, then 
	\begin{equation*}
		z(a_1)=z(b_1),\ \ \ [a_1,b_1]\subset[a,b],\ \ \ b_1-a_1=\frac{b-a}{2}.
	\end{equation*}
	By repeating the above process we get a sequence of closed intervals $\{[a_m,b_m]\}_{m=1}^{+\infty}$ satisfying 
	\begin{equation*}
		b_m-a_m=\frac{b-a}{2^m},\ \ \ z(a_m)=z(b_m),\ \ \ m=1,2,3,\cdots.
	\end{equation*}
	Therefore, we can find an unique $\xi$ such that 
	\begin{equation*}
		a_m\leq\xi\leq b_m,\ m=1,2,3,\cdots\ \ \ \text{and}\ \ \ \lim_{m\rightarrow+\infty}a_m=\lim_{m\rightarrow+\infty}b_m=\xi.
	\end{equation*}
	It is also easy to see that 
	\begin{equation*}
		\lim_{m\rightarrow+\infty}\frac{z(b_m)-z(a_m)}{b_m-a_m}=z'(\xi).
	\end{equation*}
	\qed }
	
	We remind the reader that we can make $\xi\in(a,b)$. In fact, if $a=a_1=a_2$, then 
	\begin{equation*}
		z(a)=z(a_1)=z(a_2)=z(b_1)=z(b_2).
	\end{equation*}
	One can replace $[a_2,b_2]$ by $[b_1,b_2]$ to repeat the above steps. Furthermore, by the same idea, we can make 
	$\xi\in[a+\frac{1}{8}(b-a),a+\frac{3}{8}(b-a)]$. By the classical argument we see that the $\xi$ in the Lagrange mean value theorem can also belong to such an interval.
\end{remark}

Now we can prove Lemma \ref{lemeF}.

\noindent{\it{Proof of Lemma \ref{lemeF}. }} %If $e=0$ then $F_s=0$ and there is nothing to prove. 
Following the previous argument, we just need to show \eqref{eqlemeF1} in $\Omega_2$, where 
\begin{equation}\label{eqex}
	\frac{3}{4}|x|\leq|e|\leq\frac{5}{4}|x|.
\end{equation} 
By applying the fractional Laplacian and its fundamental solution (see \cite{ah,qb,s}) we can write
\begin{equation*}
	\frac{g(x+e)-g(x)}{|e|^{2s}}=\frac{c_{n,s}}{|e|^{2s}}\left(\int_{\rn}\frac{(-\De)^sg(y)}{|x+e-y|^{n-2s}}dy
	-\int_\rn\frac{(-\De)^sg(y)}{|x-y|^{n-2s}}dy\right)
\end{equation*}
and 
\begin{equation*}
	\frac{g(x)-g(x-e)}{|e|^{2s}}=\frac{c_{n,s}}{|e|^{2s}}\left(\int_{\rn}\frac{(-\De)^sg(y)}{|x-y|^{n-2s}}dy
	-\int_\rn\frac{(-\De)^sg(y)}{|x-e-y|^{n-2s}}dy\right).
\end{equation*}

Therefore,
\begin{equation*}\label{eqlemI}
	\begin{aligned}
		\Gamma:=\frac{g(x+e)-g(x)}{|e|^{2s}}-\frac{g(x)-g(x-e)}{|e|^{2s}}
		=\frac{c_{n,s}}{|e|^{2s}}\int_{\rn}(-\De)^sg(y)\Psi dy,
	\end{aligned}
\end{equation*}
where $\Psi$ is defined in \eqref{eqpsi}. It is clear that $F_s(x)=n\left(\left(f_T(x)\right)^{\frac{n-1}{n}}\right)\Gamma$. 

%\noindent\textbf{Step 2.} We control $\Gamma$ and thus establish \eqref{eqlemeF1} in $\Omega_2$ with the aid of 
% \eqref{eqlemI}-\eqref{eqlemeFJ}. 
By the regularity result (Proposition \ref{prorf}) and the decay of fractional Laplacian 
(Proposition \ref{prodf}) we get that $(-\De)^sg\in C^{\az-2s}(\rn)$ and 
\begin{equation}\label{eqlemFgd}
	|(-\Delta)^sg(x)|\leq 
	\left\{
	\begin{aligned}
		C'_1|x|^{-\min\{\be,n\}-2s},\ \ \ \be\neq n,\\
		C'_1|x|^{-n-2s}(\ln|x|),\ \ \ \be=n,
	\end{aligned}
	\right.\ \ \ \ \ \ \text{for}\ |x|>2.
\end{equation} 
%We just consider the case $\be\neq n$ since when $\be=n$ we have $|x|^{-n-2s}(\ln|x|)\leq|x|^{-n-2s+\ep}$ for $\ep>0$ small and the proof is the same. 
Define for $|x|>2$ that 
\begin{align*}
	&A_1:=\{y\in\rn:\ |y|\leq\frac{|x|}{2}\},\\
	&A_2:=\{y\in\rn:\ |x-y|\leq\frac{|x|}{2}\},\\
	&A_3:=\rn\setminus\left(A_1\cup A_2\right).
\end{align*}
Then 
\begin{equation*}
	|\Gamma|\leq \frac{C}{|e|^{2s}}\int_{A_1\cup A_2\cup A_3}|(-\Delta)^sg(y)||\Psi|dy.
\end{equation*}

\noindent{\bf Step 1.} We assume $\beta\neq n$ and show that
\begin{equation}\label{eqlemeFgamma}
	|\Gamma|\leq\left\{
	\begin{array}{ll}
		C|x|^{-\min\{\be+2s,n\}},&\ \ \ \ \ \be+2s\neq n\\
		C|x|^{-\be-2s}\ln|x|,&\ \ \ \ \ \be+2s=n
	\end{array}
    \right.
\end{equation}
for $x\in\Omega_2$.  To do it, we control the integral in $A_1$, $A_2$ and $A_3$ respectively.

For the case in $A_1$, we first write 
\begin{gather*}
	\Gamma_1:=\int_{|y|\leq\frac{|x|}{2}}(-\De)^sg(y)\left(\frac{1}{|x+e-y|^{n-2s}}-\frac{1}{|x-y|^{n-2s}}\right)dy,\\
	\Gamma_2:=\int_{|y|\leq\frac{|x|}{2}}(-\De)^sg(y)\left(\frac{1}{|x-y|^{n-2s}}-\frac{1}{|x-e-y|^{n-2s}}\right)dy.
\end{gather*}
Note that $|x+e|\leq\frac{|x|}{4}$ in $\Omega_2$ and $(-\De)^sg$ is bounded, we have, for $\rho>0$ small, 
\begin{equation*}
	\begin{aligned}
		\Gamma_1=&\left(\int_{B_{\frac{|x|}{2}}\setminus B_\rho(x+e)}+\int_{B_\rho(x+e)}\right)
		(-\De)^sg(y)\left(\frac{1}{|x+e-y|^{n-2s}}-\frac{1}{|x-y|^{n-2s}}\right)dy\\
		\leq&C_{n,s}\int_{B_{\frac{|x|}{2}}\setminus B_\rho(x+e)}|(-\De)^sg(y)|\frac{|e|}{|\xi-y|^{n+1-2s}}dy+C_{n,s}(\rho^{2s}+\rho^n),
	\end{aligned}
\end{equation*}
where $\xi$ is obtained by the Lagrange mean value theorem. In fact, $\xi$ can be chosen such that $\xi=x+\theta e$ for $\theta\in[\frac{1}{8},\frac{3}{8}]$, 
which implies $|\xi|\geq|x|-\frac{3}{8}|e|\geq\frac{17}{32}|x|$ since $|e|\leq\frac{5}{4}|x|$. 
Moreover, $|\xi-y|\geq\frac{1}{32}|x|$ for $|y|\leq\frac{|x|}{2}$. Passing to the limit $\rho\rightarrow0$ one gets 
\begin{equation*}
	|\Gamma_1|\leq\frac{C_{n,s}|e|}{|x|^{n+1-2s}}\int_{|y|\leq\frac{|x|}{2}}|(-\De)^sg(y)|dy.
\end{equation*}
Using the same argument we also get 
\begin{equation*}
	|\Gamma_2|\leq\frac{C_{n,s}|e|}{|x|^{n+1-2s}}\int_{|y|\leq\frac{|x|}{2}}|(-\De)^sg(y)|dy.
\end{equation*}
Thus, noticing that $\frac{|x|}{2}>1$, we have 
\begin{equation}\label{eqg}
\begin{aligned}
	\int_{A_1}|(-\Delta)^sg(y)||\Psi|dy
%	&\leq\frac{C}{|x|^{n+2-2s}}\int_{B_{\frac{|x|}{2}}}|(-\Delta)^sg(y)|dy\\
	&\leq\frac{C_{n,s}|e|}{|x|^{n+1-2s}}\left(\int_{B_{1}}+\int_{B_{\frac{|x|}{2}}\setminus B_{1}}\right)|(-\Delta)^sg(y)|dy.
\end{aligned}
\end{equation}
%Since $(-\Delta)^sg$ is continuous and satisfies \eqref{eqlemFgd}, we get that $(-\Delta)^sg$ is bounded and thus   
The boundedness of $(-\Delta)^sg$ implies
\begin{equation*}
	\int_{B_{1}}|(-\Delta)^sg(y)|dy\leq C'
\end{equation*}
for some $C'>0$. On the other hand, equation \eqref{eqlemFgd} yields that 
\begin{align*}
	\int_{B_{\frac{|x|}{2}}\setminus B_{1}}|(-\Delta)^sg(y)|dy&\leq C\int_{B_{\frac{|x|}{2}}\setminus B_{1}}|y|^{-\min\{\be,n\}-2s}dy\\
	&\leq
	\left\{
	\begin{array}{ll}
		C|x|^{n-\min\{\be,n\}-2s},&\ \ \ \min\{\be,n\}+2s\neq n,\\
		C\ln|x|,&\ \ \ \min\{\be,n\}+2s=n.
	\end{array}
	\right.
\end{align*}
It follows that 
\begin{align*}
	\int_{B_{\frac{|x|}{2}}}|(-\Delta)^sg(y)|dy
	&\leq\left\{
	\begin{array}{ll}
		C(|x|^{n-\min\{\be,n\}-2s}+1),&\ \ \ \be+2s\neq n,\\
		C(\ln|x|+1),&\ \ \ \be+2s=n.
	\end{array}
	\right.
\end{align*}
Since 
\begin{align*}
	(|x|^{n-\min\{\be,n\}-2s}+1)\cdot{|x|^{-n-1+2s}}&\leq|x|^{-\min\{\be,n-2s\}-1},
\end{align*}  
we use \eqref{eqg} and \eqref{eqex} to obtain 
\begin{align*}
	\frac{C}{|e|^{2s}}\int_{A_1}|(-\Delta)^sg(y)||\Psi|dy
	&\leq\left\{
	\begin{array}{ll}
		C|e|^{1-2s}|x|^{-\min\{\be,n-2s\}-1},\\
		C|e|^{1-2s}|x|^{-n+2s-1}(\ln|x|),
	\end{array}
	\right.\\
	&\leq\left\{
	\begin{array}{ll}
		C|x|^{-\min\{\be+2s,n\}},&\ \ \ \be+2s\neq n,\\
		C|x|^{-\beta-2s}(\ln|x|),&\ \ \ \be+2s=n,
	\end{array}
	\right.\ \ \ \text{for}\ |x|>2.
\end{align*}

%The estimates in $A_2$ and $A_3$ are the same as in Lemma \ref{lemd1}. The lemma is proved.

In $A_2$ we have $|x-y|\leq\frac{|x|}{2}\leq|y|$. Define 
\begin{equation*}
	A_2^+:=\{y\in A_2:\ |x-y|\geq\frac{|x|}{4}\}\ \ \ \text{and}\ \ \ A_2^-:=A_2\setminus A_2^+.
\end{equation*}
Then \eqref{eqlemFgd} and \eqref{eqex} allow us to calculate in $A_2^+$ that  
\begin{align*}
	\frac{C}{|e|^{2s}}\int_{A_2^+}|(-\Delta)^sg(y)||\Psi|dy
	&\leq\frac{C|e|^{2-2s}}{|x|^{\min\{\be,n\}+2s+2}}\int_{\frac{|x|}{4}\leq|x-y|\leq\frac{|x|}{2}}\frac{1}{|x-y|^{n-2s}}dy\\
	&\leq C|e|^{2-2s}|x|^{-(\min\{\be,n\}+2s)+2s-2}\leq C|x|^{-\min\{\be+2s,n\}}.
\end{align*}
In $A_2^-$ we have $|x-y|<\frac{|x|}{4}$. Moreover, one finds $|x-y|<\frac{|e|}{2}$, %and thus $|x\pm e-y|>\frac{|e|}{2}>|x-y|$ in this case, 
which, combining with \eqref{eqlemeFJ} and \eqref{eqex}, implies  
\begin{align*}
	\frac{C}{|e|^{2s}}\int_{A_2^-}|(-\Delta)^sg(y)||\Psi|dy
	&\leq\frac{C}{|e|^{2s}|x|^{\min\{\be,n\}+2s}}\int_{|x-y|<\frac{|x|}{4}}\frac{1}{|x-y|^{n-2s}}dy\\
	&\leq C|e|^{-2s}|x|^{-(\min\{\be,n\}+2s)+2s}\leq C|x|^{-\min\{\be+2s,n\}}.
\end{align*}

%It allows us to use \eqref{eqlemI} to estimate $\Gamma$ since we have $|x-y|<\frac{|e|}{2}$ and thus $|x\pm e-y|>\frac{|e|}{2}>|x-y|$ in this case.

For $A_3$, we divide it into two parts
\begin{equation*}
	A_3^+:=\{y\in A_3:\ |x-y|\geq|y|\}\ \ \ \text{and}\ \ \ A_3^-:=A_3\setminus A_3^+.
\end{equation*}
Then we have $\frac{|x|}{2}\leq|y|\leq|x-y|$ in $A_3^+$ and hence $|(-\Delta)^sg(y)|\leq C|x|^{-\min\{\be,n\}-2s}$, which gives 
\begin{align*}
	\frac{C}{|e|^{2s}}\int_{A_3^+}|(-\Delta)^sg(y)||\Psi|dy&\leq
	\frac{C|e|^{2-2s}}{|x|^{\min\{\be,n\}+2s}}\int_{|y|\geq\frac{|x|}{2}}\frac{1}{|y|^{n+2-2s}}dy
	\leq C|x|^{-\min\{\be+2s,n\}}.
\end{align*}
In $A_3^-$ there holds $\frac{|x|}{2}\leq|x-y|\leq|y|$ and thus we have 
\begin{align*}
	\frac{C}{|e|^{2s}}\int_{A_3^-}|(-\Delta)^sg(y)||\Psi|dy&\leq
	\frac{C|e|^{2-2s}}{|x|^{\min\{\be,n\}+2s}}\int_{|x-y|\geq\frac{|x|}{2}}\frac{1}{|x-y|^{n+2-2s}}dy
	\leq C|x|^{-\min\{\be+2s,n\}}.
\end{align*}
We finally get \eqref{eqlemeFgamma} together with the above estimates. 

\noindent{\bf Step 2.} We consider $\beta=n$ and prove that 
\begin{equation}\label{eqG2}
	|\Gamma|\leq C|x|^{-\beta-2s}\ln|x|
\end{equation}
for $x\in\Omega_2$, by controlling the integral in $A_1$, $A_2$ and $A_3$ respectively as in {\bf Step 1}. 

In $A_1$, we still have 
\begin{equation*}
	\int_{B_{1}}|(-\Delta)^sg(y)|dy\leq C'
\end{equation*}
for some $C'>0$. By using \eqref{eqlemFgd} we obtain  
\begin{align*}
	\int_{B_{\frac{|x|}{2}}\setminus B_{1}}|(-\Delta)^sg(y)|dy&\leq C\int_{B_{\frac{|x|}{2}}\setminus B_{1}}|y|^{-n-2s}\ln|y|dy
	\leq C\left(|x|^{-2s}+1\right)\ln|x|
\end{align*}
and thus
\begin{align*}
	\int_{B_{\frac{|x|}{2}}}|(-\Delta)^sg(y)|dy\leq C\left(|x|^{-2s}+1\right)\ln|x|.
\end{align*}
Therefore, \eqref{eqg} and \eqref{eqex} imply 
\begin{align*}
	\frac{C}{|e|^{2s}}\int_{A_1}|(-\Delta)^sg(y)||\Psi|dy\leq C|e|^{1-2s}|x|^{-n-1}\ln|x|\leq C|x|^{-\beta-2s}\ln|x|\ \ \ \text{for}\ |x|>2.
\end{align*}

In $A_2$, we observe that $|y|\leq\frac{3}{2}|x|$, which implies $\ln|y|\leq C\ln|x|$. Then, the similar calculation as in 
{\bf Step 1} show that 
\begin{equation*}
	\frac{C}{|e|^{2s}}\int_{A_2^+}|(-\Delta)^sg(y)||\Psi|dy\leq C|e|^{2-2s}|x|^{-n-2}\ln|x|\leq C|x|^{-\beta-2s}\ln|x|
\end{equation*}
and 
\begin{equation*}
	\frac{C}{|e|^{2s}}\int_{A_2^-}|(-\Delta)^sg(y)||\Psi|dy\leq C|e|^{-2s}|x|^{-n}\ln|x|\leq C|x|^{-\beta-2s}\ln|x|.
\end{equation*}

In $A_3$, we first do as in {\bf Step 1} to get that 
\begin{align*}
	\frac{C}{|e|^{2s}}\int_{A_3^+}|(-\Delta)^sg(y)||\Psi|dy&\leq
	\frac{C|e|^{2-2s}}{|x|^{n+2s}}\int_{|y|\geq\frac{|x|}{2}}\frac{\ln|y|}{|y|^{n+2-2s}}dy\\
	&\leq C|e|^{2-2s}|x|^{-n-2}\ln|x|\leq C|x|^{-\beta-2s}\ln|x|.
\end{align*}
On the other hand, in $A_3^-$, it is easy to check that 
\begin{equation*}
	\ln|y|\leq\ln|x|+\ln|x-y|.
\end{equation*}
Thus, 
\begin{align*}
	&\frac{C}{|e|^{2s}}\int_{A_3^-}|(-\Delta)^sg(y)||\Psi|dy\\
	&\leq\frac{C|e|^{2-2s}\ln|x|}{|x|^{n+2s}}\int_{|x-y|\geq\frac{|x|}{2}}\frac{1}{|x-y|^{n+2-2s}}dy
	+\frac{C|e|^{2-2s}}{|x|^{n+2s}}\int_{|x-y|\geq\frac{|x|}{2}}\frac{\ln|x-y|}{|x-y|^{n+2-2s}}dy\\
	&\leq C|e|^{2-2s}|x|^{-n-2}\ln|x|\leq C|x|^{-\beta-2s}\ln|x|.
\end{align*}
\eqref{eqG2} follows from the above estimates and the lemma is proved.
\qed

Now we define 
\begin{equation*}
	h(x):=\int_\rn G(x,y)F_s^-(y)dy, 
\end{equation*}
where $G(x,y)$ denotes the fundamental solution of operator $\partial_i(a_{ij}(x)\partial_{j})$ and $F_s^-$ stands for the negative part of $F_s$. 
By \cite[Equation (7.9)]{lsw} we see that 
\begin{equation}\label{eqGdiv}
	C^{-1}|x-y|^{2-n}\leq G(x,y)\leq C|x-y|^{2-n},\ \ \ \ x,y\in\rn,\ x\neq y
\end{equation}
with constant $C>0$. Thus, $h(x)\geq0$ satisfies %\marginpar{{\color{blue}这个散度型方程应该在分布意义下来理解}}
\begin{equation*}
	\partial_i(a_{ij}(x)\partial_{j})h(x)=F_s^-(x),\ \ \ \ \forall x\in\rn.
\end{equation*}
%Note that $a_{ij}$ is divergence free (see \cite[Equation (2.6)]{tz}), we have 
By the Reilly formula (see \cite{r}) we have 
%\marginpar{\color{blue}{\tiny 这里方程如果能写成散度形式的, 就立马有\eqref{eqGdiv}, 对于非散度形式的方程, 在无界区域上Green函数的估计要比(3.13)差很多, 我们没法用. 文献\cite{py}中对非散度型方程有(3.13)这样的估计是因为其系数在无穷远处满足特定衰减条件, 我们这里不满足.}}
\begin{equation*}
	a_{ij}(x)\partial_{ij}h(x)=F_s^-(x),\ \ \ \ \forall x\in\rn,
\end{equation*}
which, combining with Lemma \ref{lemeqv}, implies that 
\begin{equation}\label{eqvsh}
	a_{ij}(x)\partial_{ij}(\De^s_ev+h)(x)\geq0,\ \ \ \ \forall x\in\rn.
\end{equation}

The decay of $h$ can be obtained by \cite[Lemma 3.2]{lb} and Lemma \ref{lemeF} since $G(x,y)$ satisfies the estimate \eqref{eqGdiv}. We really have
\begin{equation}\label{eqeh}
	0\leq h(x)\leq
	\left\{
	\begin{array}{ll}
		C(1+|x|)^{2-\min\{\be+2s,n\}},&\ \ \ \be+2s\neq n,\ \be\neq n,\\
		C(1+|x|)^{2-\beta-2s}\ln^2(2+|x|),&\ \ \ \be+2s=n,\\
		C(1+|x|)^{2-\beta-2s}\ln(2+|x|),&\ \ \ \be=n,
	\end{array}
	\right.
\end{equation}
In fact, the case $\be+2s\neq n$ and $\be\neq n$ is a direct consequence of \cite[Lemma 3.2]{lb} and Lemma \ref{lemeF}. When 
$\be+2s=n$ or $\be=n$, the decay is obtained by repeat the techniques (deal with $\ln$ term) in \cite[Lemma 3.2]{lb} and 
Lemma \ref{lemeF}.  Moreover we have the following estimate.

\begin{lemma}\label{lemvh}
	Let $s\in(0,\frac{\az}{2})$, $v$ and $h$ be as above, then 
	\begin{equation*}
		\De^s_ev(x)-h(x)\leq|e|^{2-2s},\ \ \ \forall x\in\rn,\ e\in E\setminus\{0\}.
	\end{equation*}
\end{lemma}

\noindent{\it Proof. } %If $e=0$, the inequality is trivial since $\De^s_ev(x)=0$. We always assume that $e\neq0$. 
%Note that the set $E$ consists of all the period of function $(f_p)_T$, we can find a constant $c>0$ such that 
%\begin{equation*}
%	\inf_{e\in E\setminus\{0\}}|e|\geq c>0,
%\end{equation*}
%since that matrix $T$ in Proposition \ref{pro1} satisfies $\det T=1$. 
By Remark \ref{remE} we see that 
\begin{equation*}
	\inf_{e\in E\setminus\{0\}}|e|\geq c>0.
\end{equation*}
Define for $e\in E\setminus\{0\}$ that 
\begin{equation*}
	\De^1_ev(x):=\frac{v(x+e)+v(x-e)-2v(x)}{|e|^2}.
\end{equation*}
It is easy to check that $\De^s_ev(x)=\De^1_ev(x)|e|^{2-2s}$. Therefore, by \eqref{eqvsh} we see that 
\begin{equation*}
	a_{ij}(x)\partial_{ij}(\De^1_ev+|e|^{2s-2}h)(x)\geq0,\ \ \ \ \forall x\in\rn.
\end{equation*}
Now we can do as in the proof of \cite[Equation 2.13]{tz} to get that 
\begin{equation*}
	\left(\De^s_ev(x)-h(x)\right)|e|^{2s-2}=\De^1_ev(x)-h(x)|e|^{2s-2}\leq1,
\end{equation*}
which gives the lemma.
\qed

\begin{remark}\label{remhss}
	We remind the reader that $s\in(0,\frac{\az}{2})$, $\be>2$ and $n\geq3$, thus \eqref{eqeh} gives 
	\begin{equation*}
		|h(x)|\leq C|x|^{-2s-\sigma}
	\end{equation*}
	for some $\sigma>0$.
\end{remark}

%{\marginpar{若$n=2$, 由\eqref{eqeh}可能得不到$h$的衰减(第一种情形($\be+2s\neq n,\be\neq n$)里).}}
\section{The proof of Theorem \ref{thm1}}
Now we can establish Theorem \ref{thm1}.

\noindent{\it Proof of Theorem \ref{thm1}. } 
%The proof is separated into three steps.
%
%\noindent{\bf Step 1.}

By \cite[Theorem 0.2]{cl2} (see also \cite{l}), there exists $v_p\in C^{2,\az}(\rn)$ such that 
\begin{equation}\label{eqvp}
    \left\{
    \begin{array}{ll}
    	\det(I+D^2v_p)=(f_p)_T,\ \ \ \text{in}\ \rn,\\
    	I+D^2v_p>0.
    \end{array}
    \right.
\end{equation}
Here, the function $v_p$ is periodic and has the same period as that of $(f_p)_T$. Let 
\begin{equation*}
	w(x):=v(x)-\frac{1}{2}|x|^2-v_p(x).
\end{equation*}
Applying Lemma \ref{lemvh} and the facts that $\De^s_e\left(\frac{1}{2}|x|^2\right)=|e|^{2-2s}$ and $\De^s_ev_p(x)=0$, one has
\begin{equation*}
	\De^s_ew(x)-h(x)=\De^s_ev(x)-|e|^{2-2s}-h(x)\leq0,\ \ \ \ \forall x\in\rn,\ e\in E.
\end{equation*}
Combined with the equations \eqref{eqv} and \eqref{eqvp} we have
\begin{equation*}
	\det(D^2v(x))-\det(D^2(\frac{1}{2}|x|^2+v_p(x)))=f_T(x)-(f_p)_T(x),\ \ \ \ \forall x\in\rn,
\end{equation*}
from which we get the equation for $w$:
\begin{equation*}
	\tilde{a}_{ij}(x)\pat_{ij}w(x)=f_T(x)-(f_p)_T(x),\ \ \ \ \forall x\in\rn,
\end{equation*}
where 
\begin{equation*}
	\tilde{a}_{ij}(x):=\int_0^1cof_{ij}\left[tD^2v(x)+(1-t)D^2\left(\frac{1}{2}|x|^2+v_p(x)\right)\right]dt.
\end{equation*}
In fact, we can use Proposition \ref{pro2} and \eqref{eqvp} to get that $\tilde{a}_{ij}(x)$ is uniformly elliptic. On the other hand, we can 
use the fundamental solution to construct (as that of $h$) a function $h_1$ which satisfies 
\begin{equation*}
	\tilde{a}_{ij}(x)\pat_{ij}h_1(x)=(f_p)_T(x)-f_T(x),\ \ \ \ \forall x\in\rn.
\end{equation*}
Therefore, 
\begin{equation}\label{eqwh1}
	\tilde{a}_{ij}(x)\pat_{ij}\left(w+h_1\right)(x)=0,\ \ \ \ \forall x\in\rn.
\end{equation}
By \cite[Lemma 3.2]{lb} and assumption (H) we conclude
\begin{equation*}
	|h_1(x)|\leq 
	\left\{
	\begin{array}{ll}
		C(1+|x|)^{2-\min\{\be,n\}},&\ \ \ \ \be\neq n,\\
		C(1+|x|)^{2-n}\ln(2+|x|),&\ \ \ \ \be=n.
	\end{array}
	\right. 
\end{equation*}

%\noindent{\bf Step 2.}
%We prove in this step that there exists a $b\in\rn$ such that 
%\begin{equation*}
%	|w(x)-b\cdot x|\leq C+o(|x|)\ \ \ \ \forall\ x\in\rn,
%\end{equation*}
%where $C>0$ depends only on . 

Define for $r>0$ 
\begin{equation*}
	M_r:=\sup_{B_r}(w+h_1).
\end{equation*}
We claim that, for $r$ large, there exists $M'>0$ independent of $r$ such that
\begin{equation}\label{eqtm}
	M_r\leq2M_{\frac{r}{2}}+M'.
\end{equation}

Applying Lemma \ref{lemeqv} and the definition of $w$ we see that 
\begin{equation*}
	a_{ij}(x)\pat_{ij}(\De^2_ew(x))\geq F_s(x),\ \ \ \ \forall x\in\rn.
\end{equation*}
By the construction of $h$ we have 
\begin{equation*}
	a_{ij}(x)\partial_{ij}(\De^2_ew+h)(x)\geq0,\ \ \ \ \forall x\in\rn.
\end{equation*}

Let $x_0$ be the point at which $\De^2_ew+h$ attains its maximum on $\overline{B_r}$ and $e=\frac{x_0}{2}+q$ with $|q|\leq C_q$ to be chosen so that $e\in E$ 
and $|e|\geq\frac{|x_0|}{2}$. It is clear that $x_0\in\pat B_r$. %We remind the reader that $C_q>0$ depends on $f_p$ and the linear transform $T$. 
Note that $\beta>2$, and, by Lemma \ref{lemvh}, $\De^s_ew(x)-h(x)\leq0$. Therefore, for $r$ large, 
\begin{equation*}
	\De^s_ew(x_0-e)\leq h(x_0-e)\leq C|x_0|^{-2s-\sigma}
\end{equation*} 
for some $\sigma>0$ small, which implies 
\begin{equation}\label{eqwconcave}
	w(x_0)\leq2w(x_0-e)-w(x_0-2e)+C|x_0|^{-\sigma}.
\end{equation}
By the definition of $E$, we can take $C_q=\sum_{i=1}^n|p_i|$. As a result $|w(x_0-2e)|\leq C$ since $x_0-2e\in B_{2C_q}$. Thus
\begin{equation*}
	w(x_0)\leq2\sup_{B_{\frac{r}{2}}}w+M'
\end{equation*}
and \eqref{eqtm} follows.

Next we prove by contradiction that 
\begin{equation}\label{eqtmr}
	M_r\leq \overline{M}r
\end{equation}
for some $\overline{M}>0$ independent of $r$.

Assume $\frac{M_r}{r}\rightarrow\infty$ as $r\rightarrow\infty$ and define 
\begin{equation*}
	g_r(y):=\frac{w(ry)+h_1(ry)}{M_r},\ \ \ \ \forall y\in\overline{B_1}.
\end{equation*}
Then $g_r(0)\rightarrow0$ as $r\rightarrow\infty$ since $M_r\rightarrow\infty$. We also get from \eqref{eqtm} and the definition of $M_r$ that 
\begin{equation}\label{eqglim}
	\begin{aligned}
		\max_{B_{\frac{1}{2}}}g_r\geq\frac{1}{4},&\ \ \ \ \text{for $r$ large},\\ 
		\max_{B_{1}}g_r\rightarrow1,&\ \ \ \ \text{as}\ r\rightarrow\infty.
	\end{aligned}
\end{equation}
%We observe that $1-g_r\geq0$ and $\tilde{a}_{ij}(r\cdot)\pat_{ij}(1-g_r)=0$ in $B_1$, which allows us to use Harnack inequality to get 
%\begin{equation*}
%	\max_K(1-g_r)\leq C_K\min_K(1-g_r)
%\end{equation*} 
%for any compact $K\subset B_1$. 
By \eqref{eqwh1} we find $\tilde{a}_{ij}(r\cdot)\pat_{ij}g_r=0$ in $B_1$. Applying the standard Schauder estimate we get that, 
for any compact $K\subset B_1$, $\|g_r\|_{C^\gamma(K)}$ is uniformly bounded (on $r$) for some $\gamma>0$, which yields 
that $g_r\rightarrow g$ in $C^\gamma$ norm in $K$ for a function $g$. 
%Therefore we have 
%\begin{equation*}
%	\max_K(1-g)\leq C_K\min_K(1-g).
%\end{equation*} 

On the other hand, clearly, \eqref{eqwconcave} holds for $|e|\leq\frac{3r}{4}$, which gives that $g$ is concave in $B_1$ since %$\sigma>0$ and 
the perturbation term goes to $0$ as $r\rightarrow\infty$. Let 
$l$ be a linear function touching $g$ from above at the origin. Then $l-g$ is convex and $(l-g)(0)=0$ is the minimum. By the $C^\gamma$ convergence 
we can find $l_r\rightarrow l$ such that $l_r-g_r\geq0$ in $B_{\frac{3}{4}}$ and $(l_r-g_r)(0)\rightarrow0$. 
Since $\tilde{a}_{ij}(r\cdot)\pat_{ij}(l_r-g_r)=0$, the Harnack inequality gives 
\begin{equation*}
	\max_K(l_r-g_r)\leq C_K\min_K(l_r-g_r)
\end{equation*}
for compact $K\subset B_{\frac{3}{4}}$. Therefore, combined with the definition of $g_r$ and $g$, we see that, for compact $K\subset B_1$, 
\begin{equation*}
	\max_K(l-g)\leq C_K\min_K(l-g),
\end{equation*}
which implies $l-g=0$ in $B_1$.

By \eqref{eqglim} we have $\max_{B_1}g\geq\frac{3}{4}$, thus the linear function $l=a\cdot x$ for $|a|\geq\frac{3}{4}$, where we have using the fact that 
$g(0)=l(0)=0$. It follows 
\begin{equation}\label{eqwro}
	\frac{w(ry)}{M_r}-a\cdot y\rightarrow0\ \ \ \ \text{as}\ r\rightarrow\infty\ \ \ \ \text{for }y\in B_{\frac{3}{4}},
\end{equation}
from which we obtain for $r$ large that 
\begin{equation}\label{eqcon}
	w(2^Ne)>\frac{M_r}{2r}2^N|e|,
\end{equation}
where $e\in E$ satisfies $a\cdot e>\frac{2}{3}$ and integer $N$ is chosen so that $2^N|e|<\frac{3}{4}r\leq2^{N+1}|e|$. Let $\bar{l}$ be a linear function 
that agrees with $w$ at $0$ and $e$. Then the function $\overline{w}:=w-\bar{l}$ satisfies $\overline{w}(0)=\overline{w}(e)=0$ and 
\begin{equation*}
	\overline{w}(2^Ne)>\frac{M_r}{4r}2^N|e|
\end{equation*}
for $r$ large. Since, by Remark \ref{remhss}, 
\begin{equation*}
	\De^s_e\overline{w}(x)=\De^s_ew(x)\leq C|x|^{-2s-\sigma},\ \ \ \ \forall |x|\geq1
\end{equation*} 
for some $\sigma>0$ small, we have 
\begin{equation*}
	\overline{w}(2e)\leq2\overline{w}(e)-\overline{w}(0)+C|e|^{2s}|e|^{-2s-\sigma}\leq C|e|^{-\sigma}.
\end{equation*}
Similarly, 
\begin{align*}
	\overline{w}(4e)&\leq2\overline{w}(2e)-\overline{w}(0)+C|2e|^{2s}|2e|^{-2s-\sigma}\\
	&\leq 2C|e|^{-\sigma}+C|2e|^{-\sigma}.
\end{align*}
By induction, 
\begin{align*}
	\overline{w}(2^Ne)&\leq C|e|^{-\sigma}\left(2^{N-1}+2^{N-2-\sigma}+2^{N-3-2\sigma}+\cdots+2^{-\sigma(N-1)}\right)\\
	&\leq 2^{N-1}C|e|^{-\sigma}\frac{1-2^{-N(1+\sigma)}}{1-2^{-(1+\sigma)}}\leq C2^N|e|,
\end{align*}
a contradiction to \eqref{eqcon}. Thus \eqref{eqtmr} is proved.

Now we rewrite \eqref{eqwro} as 
\begin{equation}\label{eqwo}
	\left|w(x)-\frac{M_r}{r}a\cdot x\right|=o(1)M_r,\ \ \ \ \ \text{for }x\in B_{\frac{3}{4}r}.
\end{equation}
Let 
\begin{equation*}
	b:=\lim_{r\rightarrow\infty}\frac{M_r}{r}a\ \ \ \ \text{and }\ \ \ w_1(x):=w(x)-b\cdot x. 
\end{equation*}
Then we get from \eqref{eqwo} that $|w_1(x)|=o(|x|)$ for $|x|$ large, which gives 
\begin{equation*}
	|w_1(x)|\leq C+o(|x|),\ \ \ \ \forall x\in\rn.
\end{equation*}
In other words, $w_1$ is bounded in $\rn$ and so is $w_1+h_1$. We also have 
\begin{equation*}
	\tilde{a}_{ij}(x)\pat_{ij}\left(w+h_1\right)(x)=0,\ \ \ \ \forall x\in\rn.
\end{equation*}
Thus, by the Harnack inequality we see that $w_1+h_1\equiv C$. Theorem \ref{thm1} is proved. 
\qed
%\noindent{\bf Step 3.}

\subsection*{Acknowledgments}
The authors thank the Referees for their very valuable comments.

This work was supported by %the National Key Research and Development Program of China (2020YFA0712904).
the National Natural Science Foundation of China (12371200).

\noindent\textbf{Data availability} 

\noindent Data sharing not applicable to this article as no datasets were generated or analysed for this
manuscript.

\noindent\textbf{Declarations}

\noindent{\bf Conflict of interest} The authors have no relevant financial or non-financial interests to disclose.

%\appendix
%\renewcommand{\appendixname}{Appendix~\Alph{section}}
%\section*{Appendix: An asymptotic result}
%\addcontentsline{toc}{section}{Appendix: An asympotic result}
%\renewcommand{\thesection}{\Alph{section}}
%\stepcounter{section}
%%\the\value{section}
%\newtheorem{atheorem}{Theorem}[section]

%LMAM, School of Mathematical Sciences,  
% Beijing Normal University, Beijing, 100875,
% P. R. China
%
%Shuai Qi,\quad
%E-mail address: qshuai@bnu.edu.cn


\begin{thebibliography}{99}
	
\bibitem{ah}{D. Adams, L. Hedberg, Function spaces and potential theory. Grundlehren der mathematischen 
	Wissenschaften [Fundamental Principles of Mathematical Sciences], 314. Springer-Verlag, Berlin, 1996. xii+366 pp.}

\bibitem{blz}{J. Bao, H. Li, L. Zhang, Monge-Amp\`{e}re equation on exterior domains. 
	Calc. Var. Partial Differential Equations, 52(2015), 39–63.}

%\bibitem{bxz}{J. Bao, J. Xiong, Z. Zhou, Existence of entire solutions of Monge-Ampère equations with prescribed asymptotic behavior. 
%	Calc. Var. Partial Differential Equations, 58(2019), 12 pp.}

\bibitem{b}{C. Bucur, Some observations on the Green function for the ball in the fractional Laplace framework. 
	Commun. Pure Appl. Anal., 15(2016), 657–699. }

\bibitem{ca1}{L. Caffarelli, A localization property of viscosity solutions to the Monge-Amp\`{e}re equation and their 
	strict convexity. Ann. of Math. (2), 131(1990), no. 1, 129–134.}

\bibitem{c}{L. Caffarelli, Interior $W^{2,p}$ estimates for solutions of the Monge-Amp\`{e}re equation. 
	Ann. of Math. (2), 131(1990), 135–150.}

\bibitem{ca}{L. Caffarelli, Topics in PDEs: The Monge-Amp\'{e}re equation. 
	Graduate course. Courant Institute, New York University, 1995.}

\bibitem{cl}{L. Caffarelli, Y. Li, An extension to a theorem of J\"{o}rgens, Calabi, and Pogorelov. 
	Comm. Pure Appl. Math., 56(2003), 549–583.}

\bibitem{cl2}{L. Caffarelli, Y. Li, A Liouville theorem for solutions of the Monge-Amp\`{e}re equation with periodic data. 
	Ann. Inst. H. Poincar\'{e} C Anal. Non Lin\'{e}aire, 21(2004), 97–120. }

\bibitem{ce}{E. Calabi, Improper affine hyperspheres of convex type and a generalization of a theorem by K. J\"{o}rgens. 
	Michigan Math. J., 5(1958), 105–126. }

%\bibitem{clm}{W. Chen, Y. Li, P. Ma, The fractional Laplacian. World Scientific Publishing Co. Pte. Ltd., Hackensack, NJ, 
%	[2020], ©2020. x+331 pp.}

%\bibitem{cy}{S. Cheng, S. Yau, Complete affine hypersurfaces. I. The completeness of affine metrics. 
%	Comm. Pure Appl. Math., 39(1986), 839–866.}
%
%\bibitem{d}{P. Delano\"{e}, Partial decay on simple manifolds. Ann. Global Anal. Geom., 10(1992), 3–61.}

\bibitem{dpv}{E. Di Nezza, G. Palatucci, E. Valdinoci, Hitchhiker's guide to the fractional Sobolev spaces. 
	Bull. Sci. Math., 136(2012), 521–573.}

%\bibitem{eg}{L. Evans, R. Gariepy, Measure Theory and Fine Properties of Functions, revised ed., in: Textbooks
%	in Mathematics, CRC Press, Boca Raton, FL, 2015, p. xiv+299.}

%\bibitem{fmm1}{L. Ferrer, A. Martínez, F. Mil\'{a}n, An extension of a theorem by K. J\"{o}rgens and a maximum principle at 
%	infinity for parabolic affine spheres. Math. Z., 230(1999), 471–486.}
%
%\bibitem{fmm2}{L. Ferrer, A. Martínez, F. Mil\'{a}n, The space of parabolic affine spheres with fixed compact boundary. 
%	Monatsh. Math., 130(2000), 19–27.}

\bibitem{fjm}{A. Figalli, Y. Jhaveri, C. Mooney, Nonlinear bounds in H\"{o}lder spaces for the Monge-Amp\`{e}re equation. 
	J. Funct. Anal., 270(2016), 3808–3827. }

%\bibitem{gs}{D. Gilbarg, J. Serrin, On isolated singularities of solutions of second order elliptic differential equations. 
%	J. Analyse Math., 4(1955/56), 309–340. }

%\bibitem{jw}{H. Jian, X. Wang, Continuity estimates for the Monge-Amp\`{e}re equation. 
%	SIAM J. Math. Anal., 39(2007), 608–626.}
\bibitem{jtx}{ T. Jin, X. Tu, J. Xiong, Regularity and classification of the free boundary for a Monge-Amp\`{e}re obstacle problem. 2025, arXiv:2504.21253.}


\bibitem{j}{K. J\"{o}rgens, \"{U}ber die L\"{o}sungen der Differentialgleichung $rt-s^2=1$. 
	(German) Math. Ann., 127(1954), 130–134. }

%\bibitem{jx}{J. Jost, Y. Xin, Some aspects of the global geometry of entire space-like submanifolds. 
%	Results Math., 40(2001), 233–245.}
	
\bibitem{l}{Y. Li, Some existence results for fully nonlinear elliptic equations of Monge-Amp\`{e}re type. 
    Comm. Pure Appl. Math., 43(1990), 233–271.}

\bibitem{ll}{Y. Li, S. Lu, Monge-Amp\`{e}re equation with bounded periodic data. 
	Anal. Theory Appl., 38(2022), 128–147.}

\bibitem{lb}{Z. Liu, J. Bao, Asymptotic expansion at infinity of solutions of Monge-Amp\`{e}re type equations. 
	Nonlinear Anal., 212(2021), 17 pp.}
	
\bibitem{lsw}{W. Littman, G. Stampacchia, H. Weinberger, Regular points for elliptic equations with discontinuous coefficients. 
    Ann. Scuola Norm. Sup. Pisa (3), 17(1963), 43–77.}

%\bibitem{pp}{A. Petrosyan, C. Pop, Optimal regularity of solutions to the obstacle problem for the
%	fractional Laplacian with drift. J. Funct. Anal., 268(2015), 417–472.}

\bibitem{p}{A. Pogorelov,  On the improper convex affine hyperspheres. Geometriae Dedicata, 1(1972), 33–46.}

%\bibitem{py}{Y. Pinchover, On the equivalence of Green functions of second order 
%	elliptic equations in $\rn$. Differential Integral Equations, 5(1992), 481–493.}

\bibitem{qb}{S. Qi, J. Bao, Asymptotic expansion at infinity of solutions to Monge-Amp\`{e}re equation with $C^\az$ right term. 
	2025, arXiv:2501.16667.}
	
\bibitem{r}{R. Reilly, On the Hessian of a function and the curvatures of its graph, Michigan Math. J., 20(1973), 373–383.}

\bibitem{sa}{H. Samelson, On Rolle's Theorem, The American Mathematical Monthly, 86(1979), 486–486.}

\bibitem{s}{L. Silvestre,  Regularity of the obstacle problem for a fractional power of the Laplace operator. 
	Comm. Pure Appl. Math., 60(2007), 67–112.}

%\bibitem{ss}{H. Sukjung, K. Seick, Green's function for second order elliptic equations in non-divergence form. 
%	Potential Analysis, 52(2020), 27–39.}

%\bibitem{tx}{J. Tan, J. Xiong, A Harnack inequality for fractional Laplace equations with lower order terms. 
%	Discrete Contin. Dyn. Syst., 31(2011), 975-983.}

\bibitem{tz}{E. Teixeira, L. Zhang, Global Monge-Amp\`{e}re equation with asymptotically periodic data. 
	Indiana Univ. Math. J., 65(2016), 399–422.}

\bibitem{y}{M. Yan, Extension of convex function. J. Convex Anal., 21(2014), 965–987.}
%\bibitem{aus}{D. Apushkinskaya, N. Ural'tseva, H. Shahgholian, On the Lipschitz property of the free boundary in a 
%	parabolic problem with an obstacle. Algebra i Analiz, 15 (2003), 78–103; 
%	translation in St. Petersburg Math. J., 15 (2004), 375–391.}
%

%\bibitem{c}{L. Caffarelli, Uniform Lipschitz regularity of a singular perturbation problem. 
%	Differential Integral Equations, 8 (1995), 1585–1590.}
%	


\end{thebibliography}
\end{document}